\documentclass[11pt]{article}
\usepackage{graphicx} 
\usepackage[utf8]{inputenc}
\usepackage[T1]{fontenc}
\usepackage[french]{babel}

\usepackage{newpxtext} 
\usepackage{microtype} 
\usepackage{amsmath,amsfonts,amssymb,amsthm}
\usepackage{enumitem} 
\usepackage{tikz}
\usetikzlibrary{calc,patterns,angles,quotes}

\usepackage{include}

\title{Finitude du groupe de Tate-Shafarevich pour les groupes de type multiplicatif constants sur des corps des fonctions}
\author{Melvyn El Kamel-{}-Meyrigne}
\date{}

\begin{document}
\maketitle

\textbf{Résumé} Soient $k_0$ un corps de nombres, $K$ une extension finie de $k_0(\!(x_1,...,x_n)\!)$ et soit $R$ la clôture intégrale de $k_0[[x_1,...,x_n]]$ dans $K$. Soit $G$ un groupe de type multiplicatif défini sur $k_0$. On définit le groupe de Tate-Shafarevich $\Sha^1(K,G)$ par rapport aux points de codimension $1$ de $Spec(R)$. On établit la finitude du groupe  $\Sha^1(K,G)$ lorsque $G$ provient d'un groupe de type multiplicatif $G_{k_0}$ défini sur $k_0$ sous réserve que deux hypothèses techniques soient satisfaites. On montre que ce groupe est trivial lorsque l'anneau des entiers $R$ est régulier.\\

\textbf{Abstract} Let $k_0$ be a number field, let $K$ be a finite extension of $k_0(\!(x_1,...,x_n)\!)$ and let $R$ be the integral closure of $k_0[[x_1,...,x_n]]$ in $K$. Consider a group of multiplicative type $G$ defined over $K$. We define the Tate-Shafarevich group $\Sha^1(K,G)$ with respect to the points of codimension $1$ of $Spec(R)$. We show the finiteness of the group $\Sha^1(K,G)$ when $G$ comes from a group of multiplicative type $G_{k_0}$ defined over $k_0$ provided that two technical conditions are satisfied. We then prove that this group is trivial when the ring of integers $R$ is regular.\\

\textbf{Mots clés :} Groupes de Tate-Shafarevich, Anneaux complets, Tores algébriques, Cohomologie galoisienne.

\textbf{Keywords} : Tate-Shafarevich
groups, Complete ring, Algebraic torus, Galois cohomology.

\section{Introduction}

Afin de comprendre la structure des points rationnels d'une variété définie sur un corps de nombres $k$, une approche particulièrement féconde consiste à tenter d'établir l'existence de points locaux, c'est-à-dire sur les complétés de $k$ (qui sont des corps dont l'arithmétique est mieux comprise), puis d'en déduire l'existence d'un point global sur $k$. Si l'on se donne une famille $\sF$ de variétés $Z$ définies sur $k$ et que l'on note $\Omega_k$ l'ensemble des places de $k$, on dit que $\sF$ satisfait le principe local-global si
$$
\forall Z \in \sF, \prod_{v \in \Omega_k}Z(k_v) \neq \emptyset \Ra Z(k) \neq \emptyset.
$$
L'exemple le plus notoire de variétés satisfaisant le principe local-global est donné par les quadriques définies sur un corps de nombres d'après le théorème de Hasse-Minkowski. Cependant, ce principe est loin d'être satisfait de manière générale ce qui nous amène à étudier les possibles obstructions à sa véracité.\\

Soit $G$ un groupe algébrique commutatif lisse sur $k$. Le groupe $H^1(k,G)$ classifie les $G$-espaces principaux homogènes à $k$-isomorphisme près. De plus, la classe d'un $G$-torseur est nulle dans $H^1(k,G)$ si et seulement si il possède un point rationnel. Ainsi, le groupe de Tate-Shafarevich
$$
\Sha^1(k,G) \coloneqq Ker(H^1(k,G) \ra \prod_{v \in \Omega_k}H^1(k_v,G))
$$
mesure exactement l'obstruction au principe local-global pour les $G$-torseurs.\\

De nombreux résultats sur ces groupes ont déjà été obtenus, l'un des plus célèbres étant le théorème de dualité de Poitou-Tate en 1967 pour les groupes finis étales et pour les tores. Plus récemment, plusieurs théorèmes de dualité similaires à celui de Poitou-Tate ont pu être démontrés pour des tores définis d'autres corps tels que le corps de fonctions $K =k_0(C)$ d'une $k_0$-courbe $C$ lorsque $k_0$ est un corps $p$-adique \cite{HarSza2}, une extension finie de $\C(\!(t)\!)$ \cite{Colliot_Th_l_ne_2015} ou plus généralement un corps $d$-local \cite{Izquierdo2}.  Dans ce cadre, on ne considère pas toutes les places du corps $K$ mais seulement celles dites géométriques, i.e. provenant de l'ensemble $C^{(1)}$ des points de codimension $1$ de $C$.\\

Cependant, on ne dispose pas d'une compréhension aussi poussée de ces groupes pour tout type de corps. De fait, les groupes de Tate-Shafarevich d'un groupe algébrique définis sur le corps des fonctions d'une variété sur un corps de nombres sont plus élusifs. Une première avancée majeure dans l'étude de ces groupes, obtenue par M.Saïdi et A.Tamagawa dans \cite{SaïdiTamagawa}, est la finitude des groupes $\Sha^1_{C^{(1)}}(K,A) \coloneqq Ker(H^1(K,A) \ra \prod_{v \in C^{(1)}}H^1(K_v,A))$ où $K=k_0(C)$ est le corps des fonctions d'une courbe $C$ définie sur un corps $k_0$ finiment engendré sur $\Q$ et $A$ est une variété abélienne isotriviale sur $K$.
D. Harari et T. Szamuely ont ensuite établi dans \cite{HarSza} la finitude du groupe $\Sha^1_{C^{(1)}}(K,G)$ lorsque $k_0$ est un corps de nombres et $G$ un groupe de type multiplicatif constant. \\

Indépendamment, A. Rapinchuk et I. Rapinchuk ont établi de manière simultanée un résultat similaire dans leur article \cite{Rapinchuk1} en faisant appel à des techniques adéliques développées dans \cite{Rapinchuk3}. Ils étudient les groupes de Tate-Shafarevich de tores constants définis sur le corps des fonctions $K$ d'un schéma normal $X$ de type fini sur $\Z$ dont le corps de fonctions est $K$, le schéma $X$ étant de dimension quelconque. Cependant, les places de $K$ qui entrent en jeu sont celles correspondant aux points de codimension $1$ de $X$; cet ensemble de places est plus large que celui considéré dans \cite{HarSza}.
Peu de temps après, les deux auteurs sont parvenus à améliorer leurs résultats dans l'article \cite{Rapinchuk2} en combinant les techniques développées dans \cite{Rapinchuk1} et les résultats de \cite{HarSza}. Ils ont ainsi réussi à réduire l'ensemble des places considérées à celles provenant des points de codimension $1$ de la fibre générique de $X$; ce faisant, ils parviennent à généraliser le résultat de Harari et Szamuely aux corps des fonctions de variétés de dimension arbitraire.\\

Introduisons à présent le cadre de cet article. Soit $k_0$ un corps de nombres et soit $K$ une extension finie de $k_0(\!(x_1,x_2,...,x_n)\!)$. Notons $R$ la clôture intégrale de $k_0[[x_1,x_2,...,x_n]]$ dans $K$ et $k$ le corps résiduel de $R$. Le but de cet article est d'obtenir un résultat de finitude sur le groupe de Tate-Shafarevich $\Sha^1_{R^{(1)}}(K,G)$ d'un groupe de type multiplicatif $G$ défini sur $k_0$ par rapport aux valuations induites par les points de codimension $1$ de $R$ en combinant des outils employés dans \cite{HarSza} et \cite{CTOP}.\\

Le corps $K$ que l'on étudie est le corps des fonctions du schéma $Spec(R)$; ce dernier joue donc le rôle de la courbe sur $k_0$ considérée dans \cite{HarSza}, la différence fondamentale étant que  $Spec(R)$ n'est ni de type fini sur $k_0$, ni un schéma régulier. Afin de pallier les difficultés introduites par les éventuelles singularités de $R$, on travaille avec une désingularisation $\CX \ra Spec(R)$ dont la fibre singulière est un diviseur à croisement normaux strict.\\

Le point saillant de la preuve consiste à comparer le groupe de Tate-Shafarevich $\Sha^1(K,G)$ et le groupe
$$
\Sha_{n}^1(k,G) \coloneqq Ker(H^1(k,G) \ra \prod_{c \in \CX^{(n)}}H^1(k(c),G))
$$
où $\CX^{(n)}$ désigne l'ensemble des points fermés de $\CX$, la structure de ce dernier groupe étant mieux comprise.
Deux hypothèses, l'une de nature arithmétique et l'autre de nature géométrique, sont nécessaires au bon déroulement de la preuve. La première permet de garantir la finitude du groupe $\Sha_{n}^1(k,G)$.

\begin{hypoe}\label{hyp1}
Soient $Z_1,...,Z_l$ les composantes irréductibles de la fibre spéciale $\CX_k$ munie de sa structure réduite et posons $k_i \coloneqq k(Z_i) \cap \overline{k}$. On suppose que le noyau du morphisme
$$
H^1(k, G) \ra \prod_{i=1}^l H^1(k_i, G)
$$
est fini.
\end{hypoe}

La seconde hypothèse permet de garantir que la cohomologie de la fibre spéciale ne diffère pas trop de celle de ses composantes irréductibles.
\begin{hypoe}\label{hyp2}
Soient $Z_1,...,Z_l$ les composantes irréductibles de la fibre spéciale $\CX_k$ munie de sa structure réduite. On suppose que le noyau du morphisme $\Pic(Z) \ra \prod_{i=1}^l\Pic(Z_i)$ est fini.
\end{hypoe}
En dimension $2$, ces hypothèses sont automatiquement satisfaites lorsque la singularité de $R$ rationnelle. 

On peut à présent énoncer le résultat principal de l'article :

\begin{thme}
     Supposons que les hypothèses \ref{hyp1} et \ref{hyp2} sont satisfaites. Alors, le groupe $\Sha^1(K,G)$ est fini.
\end{thme}

L'hypothèse de complétude sur $R$ est essentielle : c'est elle qui permet de ramener l'étude de la cohomologie de la désingularisation à celle de la cohomologie des épaississements de sa fibre spéciale; cette approche a notamment pu être utilisée pour étudier les groupes de Brauer d'un anneau hensélien complet de dimension $2$ dans \cite{CTOP}. Dans notre cadre, la fibre spéciale est une variété sur un corps de nombres ce qui nous permet de faire appel à certains théorèmes de finitude tels que celui de Mordell-Weil et celui de la finitude des groupes des classes d'idéaux.

Une autre question naturelle est celle de la trivialité du groupe $\Sha^1(K,G)$. Le théorème \ref{thmezero} donne un premier élément de réponse.

\begin{thme}
Supposons que l'anneau $R$ soit régulier. Alors, le groupe $\Sha^1(G)$ est trivial.
\end{thme}

À la fin de l'article, on montre que l'approche utilisée dans la première partie s'adapte naturellement dans un contexte similaire : celui des corps des fonctions de variétés définies sur des extensions finies de $\Q(\!(x)\!)$, par exemple $\Q(\!(x)\!)(t_1,...,t_m).$\\

\textbf{Remerciements :} Je tiens à adresser mes remerciements les plus sincères à Diego Izquierdo pour toutes les discussions que nous avons pu avoir, son écoute permanente, sa relecture attentive et les multiples conseils qu'il m'a apportés; nul doute que ce travail n'existerait pas sans son soutien. 

\section{Notations et résultats préliminaires}
On utilisera les notations suivantes : 
\begin{enumerate}[label={$\bullet$}]
    \item $k_0$ est un corps de nombres.
    \item $k_0^s$ désigne une clôture séparable de $k_0$.
    \item $K$ est une extension finie de $k_0(\!(x_1,...,x_n)\!)$.
    \item $R$ désigne l'anneau des entiers de $K$, i.e. la clôture intégrale de $k_0[[x_1,...,x_n]]$ dans $K$. On note $\fm$ l'idéal maximal de $R$ et $k$ son corps résiduel.
    \item $R_{sing}$ désigne le fermé de $Spec(R)$ associé à l'ensemble des points singuliers de $Spec(R)$.
    \item Pour tout point $v$ de $Spec(R)$ de codimension 1, $K_v$ désigne la complétion de $K$ par rapport à la valuation $v$.
    \item $G_k$ est un groupe multiplicatif défini sur $k_0$ que l'on étend à $G \coloneqq G_{k_0} \times_{k_0} K$. On utilisera également la notation $G$ afin de désigner $G_{k_0}$ pour alléger les notations.
\end{enumerate}

De manière équivalente, on peut se donner une $k_0$-algèbre $R$ locale,
complète, normale et excellente de dimension $n$ dont le corps résiduel est un corps de nombres $k$ et considérer son corps des fractions $K$.\\

\underline{Groupes abéliens :} Soit $A$ un groupe abélien. On utilisera la notation suivante :
\begin{itemize}[label={$\bullet$}]
    \item ${}_{n}A$ : la $n$-torsion de $A$ où $n \in \N^{*}$.
    \item $\overline{A}$ : le quotient de $A$ par son sous-groupe divisible maximal.
\end{itemize}

On utilisera à multiples reprises le fait que les groupes divisibles sont les objets injectifs dans la catégorie des groupes abéliens, i.e. toute suite exacte de groupes abéliens
$$
0 \ra D \ra M \ra N \ra 0
$$
avec $D$ divisible est scindée (cf \cite{Baer}, théorème 1).\\

\underline{Schémas :} Soit $Z$ un schéma noethérien. Pour tout entier naturel $i$, notons $Z^{(i)}$ l'ensemble des points de codimension $i$ de $Z$.
Pour tout point $c \in Z$, on note $\CO^h_{Z,c}$ l'henselisé de l'anneau local $\CO_{Z,c}$. De plus, on note $K^h_c$ le corps des fractions de $\CO^h_{Z,c}$.\\
Soit $\beta_1,...,\beta_n$ un ensemble de points de $Z$. On dit que $\beta = (\beta_1,...,\beta_n)$ est un drapeau dans $Z$ si $\ol{\{\beta_{i+1}\}} \subseteq \ol{\{\beta_{i}\}}$ pour tout $1 \leq i \leq n$. De plus, on dit que $\beta$ est un drapeau régulier si les anneaux $\CO_{Z,\beta_n}/\beta_i$ sont réguliers.\\

\underline{Résolution des singularités :} Soit $\CX$ un modèle régulier de $Spec(R)$, i.e. , un schéma intègre régulier muni d'un morphisme birationnel projectif $f : \CX \ra Spec(R)$  dont la fibre singulière $\CX_{sing} \coloneqq f^{-1}(R_{sing})$ est un diviseur à croisements normaux strict et dont la fibre spéciale est également un diviseur ; un tel $\CX$ existe toujours d'après l'article \cite{temkinDesingularizationQuasiexcellentSchemes2008} qui généralise la résolution des singularités de Hironaka à des schémas noethériens quasi-excellents de caractéristique nulle. On se fixe donc un tel modèle régulier $\CX$ pour la suite.\\

On note $X^{(1)}$ l'ensemble des points de codimension $1$ de $\CX$ correspondant à des points de codimension $1$ de $R$; une telle identification est possible par normalité de $R$. Ainsi, les points de codimension $1$ de $\CX$ sont séparés en deux parties : ceux qui appartiennent à $X^{(1)}$ et des points additionnels qui correspondent à des points génériques de la fibre singulière $\CX_{sing}$ de $\CX$; seuls ceux provenant de $Spec(R)$ nous intéressent.\\

Par normalité de $R$, le théorème de factorisation de Stein implique l'identité $f_*\CO_{\CX}=\CO_R$ et les fibres de $f$ sont géométriquement connexes (cf \cite{PMIHES_1961__11__5_0}, théorème 4.3.1 et remarque 4.3.4).

Par souci pratique, on donne ici les hypothèses nécessaires au théorème afin de pouvoir le réutiliser dans d'autres contextes; on trouvera une preuve de ce résultat dans  \cite[\href{https://stacks.math.columbia.edu/tag/0AY8}{Lemma 0AY8}]{stacks-project}.

\begin{lemme}\label{lemstein}
    Soit $f : X \ra S$ un morphisme de schémas tel que
\begin{itemize}
    \item $f$ est propre.
    \item S est intègre de point générique $\xi$.
    \item S est normal.
    \item X est réduit.
    \item Chaque point générique de $X$ est envoyé sur $\xi$.
    \item $H^0(X_{\xi},\CO_X) = k(\xi).$
\end{itemize}    
Alors, $f_*\CO_X = \CO_S$ et les fibres de $f$ sont géométriquement connexes.
\end{lemme}

Commençons par établir quelques lemmes qui serviront par la suite.

\begin{lemme}\label{lemdimension}
    Soit $x \in \CX$. Alors, il existe un point fermé $\beta_n \in \CX$ tel que $\beta_n \in \ol{\{x\}}$. En outre, si l'on note $y$ l'image de $x$ par $f : \CX \ra Spec(R)$, on a que 
 $$
dim(\CO_{\CX,x}) = dim(R_{y}) - deg.tr_{k(y)}k(x),
$$
En particulier, on a que $codim(x) \leq codim(y).$  
\end{lemme}

\proof
Premièrement, comme $f$ est fermé et que $R$ local, l'ensemble fermé $f(\ol{\{x\}})$ contient le point fermé $\fm$ de $Spec(R)$; il en découle que l'ensemble $\ol{\{x\}} \cap \CX_k$ est un fermé non-vide de $\CX_k$, il contient donc un point fermé $\beta_n$ de $\CX_k$.\\
Deuxièmement, puisque $f$ est birationnel et que $R$ est un anneau universellement caténaire car noethérien, local et complet, la formule de la dimension (\cite{PMIHES_1965__24__5_0}, §5, proposition 5.6.5) fournit l'égalité voulue.
\qed
\begin{lemme}\label{lemannvaldis}
    Soient $R$ un anneau local, $S$ un anneau valuation discrète, $k_R$ et $k_S$ leurs corps résiduels respectifs et soit $\varphi : R \ra S$ un morphisme local tel que $S$ est une $R$-algèbre de type fini. Alors, $\varphi^{-1}(\{0\}) \neq \fm $.
\end{lemme}

\proof Soit $\pi_S$ une uniformisante de $S$ et notons $K_S$ le corps de fractions de $S$. Supposons par l'absurde que $\varphi^{-1}(\{0\}) = \fm $. Alors, $S$ est une $k_R$-algèbre de type fini et à fortiori, $S[T]/(\pi_ST-1) \simeq K_S$ est une $k_R$-algèbre de type finie; l'extension $K_S/k_R$ est donc finie par le lemme de Zariski ce qui est une contradiction puisque l'uniformisante $\pi_S \in K_S$ est transcendante sur $k_R$. \qed \\

\section{Relation entre les groupes $\Sha^1(K,G)$ et $\Sha^1_n(k,G)$}
Cette section a pour objectif de comparer la cohomologie des points de codimension $1$ de $Spec(R)$ et celle des points fermés de $\CX$ via le groupe $\Sha^1_n(k,G)$ puis d'étudier la structure de ce dernier groupe.
Le prochain lemme est instrumental afin de pouvoir gérer les points de codimension $1$ supplémentaires de $\CX$ provenant de la fibre singulière.

\begin{lemme}\label{lemdrapeau}
Pour tout point fermé $\alpha_n$ de $\CX$, il existe un modèle régulier $\CX'$ de $\CX$, un drapeau régulier $\beta = (\beta_1,...,\beta_n)$ dans $\CX'$ tel que $k(\alpha_n) \simeq k(\beta_n)$ et une famille $\{S_2,...,S_n\}$ d'anneaux de valuation discrète henséliens tels que 
\begin{itemize}
    \item Pour tout $2\leq i < n$, on a  $k(\beta_i) \simeq Frac(S_{i+1})$.
    \item  Pour tout $2 < i \leq n$, on a $k(S_i) \simeq k(\beta_i)$ où $k(S_i)$ désigne le corps résiduel de $S_i$.
    \item $\beta_1 \in X^{(1)}$.
\end{itemize}

\end{lemme}

La première partie de la démonstration de ce lemme reprend des arguments de la preuve du lemme 4.6 de \cite{Wittenberg}.

\proof

Rappelons que $X^{(1)}$ désigne l'ensemble des points de codimension $1$ de $\CX$ correspondant à des points de codimension $1$ de $R$. La difficulté de ce lemme réside dans la construction d'un drapeau $(\beta_1,...,\beta_n)$ où $\beta_1$ n'est pas dans la fibre singulière $\CX_{sing}$ de $\CX$.\\
La preuve est séparée en deux parties : dans un premier temps, on montre qu'il existe un drapeau régulier $(\beta_1,...,\beta_n)$ dans un modèle $\CX'$ de $\CX$ tel que chaque $\beta_i$ appartient à la fibre d'un élément de $Spec(R)$ de codimension égale à $i$. On construit ce faisant une famille d'anneaux de valuation discrète $(R_2,...,R_n)$ et on vérifie dans un second temps que leurs hensélisés satisfont les propriétés indiquées dans l'énoncé.\\

\underline{Étape 1 :} Soit $\alpha_n$ un point fermé de $\CX$.  Comme $\alpha_{n}$ est un point régulier de $\CX$, on peut, quitte à éclater $\alpha_{n}$ et à le remplacer par un point $\beta_{n}$ de même corps résiduel, supposer que $\beta_{n}$ n'est pas dans l'intersection de deux composantes irréductibles de la fibre spéciale. Le point $\beta_{n}$ est un point régulier de la fibre spéciale et la composante irréductible à laquelle il appartient est de codimension $1$.\\
Soit alors $t \in \CO_{\CX,\beta_{n}}$ tel que $Z \coloneqq (\CX_k)_{red}$ ait pour équation locale $t=0$ en $\beta_n$. Étant donné que $t$ n’est pas un diviseur de zéro dans $\CO_{\CX,\beta_n}$ et comme on a supposé que $\CO_{\CX,\beta_n} / (t) = \CO_{\CZ, \beta_n}$ est régulier, on peut, d'après le corollaire 17.1.8 et la  proposition 17.1.7 de \cite{PMIHES_1964__20__5_0}, compléter $(t)$ en un système régulier de paramètres $(t,f_2,...,f_n)$ de $\CO_{\CX,\beta_n}$ et on définit alors $\beta_{n-1} \coloneqq (f_2,...,f_n) \in \CX$ ; c'est bien un point de codimension $n-1$ par (\cite{Matsumura} chapitre 7, théorème 36).\\

On souhaite à présent montrer que $\beta_{n-1}$ est envoyé sur un point $\gamma_{n-1}$ de $Spec(R)$ tel que $codim(\gamma_{n-1}) = n-1$. Par le lemme \ref{lemdimension}, cela revient à démontrer que $\beta_{n-1} \notin \CX_k$. Pour cela, posons $R_n \coloneqq \CO_{\CX,\beta_n}/\beta_{n-1}$, c'est un anneau régulier, local de dimension $1$ : c'est donc un anneau de valuation discrète. Considérons le $k$-espace vectoriel
$$R_n \otimes_R k
 \simeq R_n \otimes_R R/\fm \simeq R_n/\fm R_n .$$

Soit $r$ la multiplicité géométrique de la composante irréductible contenant $\beta_n$, alors $ t^r  \in \fm \CO_{\CX,\beta_n}$. Cela implique que l'image de $\fm$ dans $R_n$ est non nulle et il existe donc un $l>0$ tel que $\fm R_n = (t^l)$. Il en résulte que le quotient $\simeq R_n/\fm R_n $ est alors un $k(\beta_n)$-espace vectoriel de dimension $l$. Comme $k(\beta_n)$ est une extension finie de $k$, on déduit que $R_n/\fm R_n $ est un $k$-espace vectoriel de dimension finie.\\

Comme $R_n$ est essentiellement de type fini en tant que $R$-module, il existe un ouvert affine $Spec(A)$ de $\CX$ contenant $\beta_{n}$ pris suffisamment petit pour qu'il contienne les $f_2,...,f_n$ et que $A$ soit de type fini sur $R$. Étant donné que $R$ est hensélien et que $R_n \otimes_R k$ est de dimension finie, le théorème 18.5.11, c' de \cite{PMIHES_1967__32__5_0}  appliqué à $Spec(A/f_2,...,f_n)$ montre que $R_n$ est une $R$-algèbre finie. Cela implique
que le point $\beta_{n-1}$, qui est l'image du point générique $\eta_n$ de Spec ($R_n$), est envoyé sur un point $\gamma_{n-1}$ de $Spec(R)$ de codimension $n-1$ par les lemmes \ref{lemannvaldis} et \ref{lemdimension}; c'est ce que l'on voulait.\\

La construction du point $\beta_{n-2}$ est peu ou prou identique une fois que l'on se place dans le bon cadre. Considérons le schéma  $\CX_{R_{\gamma_{n-1}}} \coloneqq \CX \times_R R_{\gamma_{n-1}}$. Remarquons que $R_{\gamma_{n-1}}$ est la limite du système direct $(R_s,f_{st})_{s \notin \gamma_{n-1}}$ où $f_{st}$ désigne le morphisme naturel $f_{st} : R_s \ra R_t$ quand $s \mid t$.

Ainsi, on a que $Spec(R_{\gamma_{n-1}}) = \prl_{s \notin \gamma_{n-1}} Spec(R_s)$ (\cite{PMIHES_1966__28__5_0}, proposition 8.2.3) . De plus, par la proposition 8.2.5 de \cite{PMIHES_1966__28__5_0}, $\CX_{R_{\gamma_{n-1}}}$ s'identifie à la limite projective des ouverts $\CX \times_R R_s$ pour $s \notin \gamma_{n-1}$. De fait, on peut voir $\CX_{R_{\gamma_{n-1}}}$ comme une partie de $\CX$ : d'une part, si l'on regarde les espaces topologiques sous-jacents, $\CX_{R_{\gamma_{n-1}}}$ n'est rien d'autre que l'ensemble des généralisations des points de $\CX_{\gamma_{n-1}}$.\\
D'autre part, pour tout $\alpha \in \CX_{R_{\gamma_{n-1}}}$, on a que $\CO_{\CX_{R_{\gamma_{n-1}}}, \alpha} \simeq \CO_{\CX, \alpha}$.
En effet, si l'on note $\varphi_s : \CX_{R_{\gamma_{n-1}}} \ra \CX_s$ le morphisme canonique pour tout $s \notin \gamma_{n-1}$, on a que
\begin{align*}
\CO_{\CX_{R_{\gamma_{n-1}}}, \alpha} & \simeq \drl_{\underset{U \text{ouvert}}{U \ni \alpha}} \CO_{\CX_{R_{\gamma_{n-1}}}}(U) \\
& \simeq \drl_{{\underset{V_s \text{ouvert affine } \CX_s}{V_s \ni f_s(\alpha)}}} \CO_{\CX_{R_{\gamma_{n-1}}}}(\varphi^{-1}_s(V_s)) \hspace{1cm} \text{par cofinalité des ouverts affines de la forme $\varphi^{-1}_s(V_s)$} \\
& \simeq \drl_{\underset{V_s \text{ouvert affine } \CX_s}{V_s \ni f_s(\alpha)}} \drl_{t \geq s} \CO_{\CX_t}(f_{ts}^{-1}(V_s)) \\
& \simeq \drl_{t \notin \gamma_{n-1}} \drl_{\underset{V_s \text{ouvert affine } \CX_s, t \geq s}{V_s \ni f_s(\alpha)}} \CO_{\CX_t}(f_{ts}^{-1}(V_s)) \\
& \simeq \drl_{t \notin \gamma_{n-1}} \CO_{\CX_t, \alpha} \\
& \simeq \CO_{\CX, \alpha}.
\end{align*}

Ainsi, par le lemme \ref{lemdimension}, $\beta_{n-1}$ correspond à un point fermé de la fibre $(\CX_{R_{\gamma_{n-1}}})_{\gamma_{n-1}}$ de $\CX_{R_{\gamma_{n-1}}}$.\\

Deux situations peuvent survenir :
\begin{itemize}
    \item Si $\gamma_{n-1}$ est un point régulier de $Spec(R)$, alors $\beta_{n-1}$ s'identifie à $\gamma_{n-1}$ et on peut simplement considérer n'importe quel système régulier de paramètres $(g_1,...,g_{n-1})$ de $\CO_{\CX,\beta_{n-1}}$ et définir chaque $\beta_i$ par $\beta_i \coloneqq (g_1,...,g_{i})$ et $R_i \coloneqq \CO_{\CX,\beta_i}/\beta_{i-1}$; on a alors automatiquement que $\beta_1 \in X^{(1)}$ .
    \item Si $\gamma_{n-1}$ est un point singulier de $Spec(R)$, on peut, comme cela a été fait pour définir $\beta_{n-1}$, se ramener au cas où $\beta_{n-1}$ est un point de $Z' \coloneqq ((\CX_{R_{\gamma_{n-1}}})_{\gamma_{n-1}})_{red}$ n'appartenant qu'à une unique composante irréductible d'équation locale $t'=0$ en $\beta_{n-1}$. On considère ensuite un système régulier de paramètres $(t',h_1,...,h_{n-2})$ et l'on pose alors $\beta_{n-2} \coloneqq (h_1,...,h_{n-2})$ et $R_{n-1} \coloneqq \CO_{\CX,\beta_{n-1}}/\beta_{n-1}$. Pour montrer que $\beta_{n-2}$ est envoyé sur un point $\gamma_{n-2}$ de $Spec(R)$ de codimension $n-2$, on procède de manière analogue à ce qui précède à une différence près : comme l'anneau $R_{\gamma_{n-1}}$ n'est pas nécessairement hensélien, on doit travailler avec son hensélisé $R_{\gamma_{n-1}}^h$.\\

Plus précisément, soit $B$ un ouvert affine de $\CX_{R_{\gamma_{n-1}}}$ contenant $\beta_{n-1}$ pris suffisamment petit pour qu'il contienne $h_1,...,h_{n-2}$ et notons $C \coloneqq B/(h_1,...,h_{n-2})$. Considérons le diagramme commutatif suivant donné par le produit fibré :

\begin{equation}
\begin{tikzcd}
    Spec(C \otimes_{R_{\gamma_{n-1}}} R_{\gamma_{n-1}}^h) & Spec(R_{\gamma_{n-1}}^h) \\
    Spec(C) & Spec(R_{\gamma_{n-1}}). \\
\arrow[from=1-1, to=1-2, "\psi"]
\arrow[from=1-1, to=2-1, "f_{n-1}"]
\arrow[from=1-2, to=2-2]
\arrow[from=2-1, to=2-2, "g"]
\end{tikzcd}
\end{equation}

Comme $\beta_{n-1}$ est un point fermé de $\CX_{R_{\gamma_{n-1}}}$, le $k(\gamma_{n-1})-$espace vectoriel $R_{n-1} \otimes_{R_{\gamma_{n-1}}} k(\gamma_{n-1})$ est de dimension finie ce qui implique que le morphisme $Spec(C) \ra Spec(R_{\gamma_{n-1}})$ est quasi-fini en $\beta_{n-1}$. En outre, l'inclusion $R_{\gamma_{n-1}} \subseteq R_{\gamma_{n-1}}^h$ étant fidèlement plate (cf \cite{PMIHES_1967__32__5_0}, théorème 18.6.6, \textit{iii}), le morphisme induit $Spec(R_{\gamma_{n-1}}^h) \ra Spec(R_{\gamma_{n-1}})$ est surjectif et il en va donc de même du morphisme $\psi : Spec(C \otimes_{R_{\gamma_{n-1}}} R_{\gamma_{n-1}}^h) \ra Spec(C)$. Soit alors $\lambda_{n-1}$ un élément de $\psi^{-1}(\beta_{n-1})$, le morphisme $Spec(C \otimes_{R_{\gamma_{n-1}}} R_{\gamma_{n-1}}^h) \ra Spec(R_{\gamma_{n-1}}^h)$ est quasi-fini au point $\lambda_{n-1}$ par changement de base ( \cite[\href{https://stacks.math.columbia.edu/tag/00PP}{Tag 00PP}]{stacks-project}).

De plus, l'image de $\lambda_{n-1}$ dans $Spec(R_{\gamma_{n-1}}^h)$ est l'idéal maximal de $R_{\gamma_{n-1}}^h $ en vertu de la commutativité du diagramme $(1)$ et de la platitude du morphisme $Spec(R_{\gamma_{n-1}}^h) \ra Spec(R_{\gamma_{n-1}})$. On peut donc utiliser le théorème 18.5.11, c' de \cite{PMIHES_1967__32__5_0} pour en déduire que la tige de $C \otimes_{R_{\gamma_{n-1}}} R_{\gamma_{n-1}}^h$ en $\lambda_{n-1}$ ,qui n'est autre que l'hensélisé de $R_{n-1}$ 
(\cite[\href{https://stacks.math.columbia.edu/tag/05WP}{Tag 05WP}]{stacks-project}) est une $R_{\gamma_{n-1}}^h$-algèbre finie.\\

En outre, on peut par platitude relever la généralisation $\beta_{n-2} \rightsquigarrow \beta_{n-1}$ en une généralisation $\lambda_{n-2} \rightsquigarrow \lambda_{n-1}$ dans $Spec(C \otimes_{R_{\gamma_{n-1}}} R_{\gamma_{n-1}}^h)$ et à fortiori dans $(R_{n-1} \otimes_{R_{\gamma_{n-1}}} R_{\gamma_{n-1}}^h)_{\lambda_{n-1}} \simeq R_{n-1}^h$. Comme $R_{n-1}$ est un anneau de valuation discrète, il en va de même de son hensélisé $R_{n-1}^h$. Ainsi, le lemme \ref{lemannvaldis} montre que l'image de $\lambda_{n-2}$ dans $Spec(R_{\gamma_{n-1}}^h)$ n'est pas le point fermé $\fm R_{\gamma_{n-1}}^h$. Par platitude, son image dans $Spec(R_{\gamma_{n-1}})$ est un point de codimension $< n-1$; le diagramme commutatif $(1)$ montre alors que $\gamma_{n-2} = f(\beta_{n-2})$ est de codimension $n-2$ par le lemme \ref{lemdimension}.
\end{itemize}

On procède de manière identique pour définir les points $\beta_{n-3},...,\beta_{1}$. Enfin, comme tous les points de codimension $1$ de $Spec(R)$ sont réguliers par normalité, on a que $\beta_1 \in X^{(1)}$.\\

\underline{Étape 2 :} On souhaite maintenant vérifier que les anneaux de valuation discrète $S_i = R_i^h = (\CO_{\CX,\beta_{i}}/\beta_{i-1})^h \simeq \CO_{\CX,\beta_{i}}^h/\beta_{i-1}\CO_{\CX,\beta_{i}}^h$ satisfont les propriétés suivantes :
\begin{itemize}
    \item Pour tout $2\leq i < n$, on a  $k(\beta_i) \simeq Frac(S_{i+1})$.
    \item  Pour tout $2 < i \leq n$, on a $k(S_i) \simeq k(\beta_i)$ où $k(S_i)$ désigne le corps résiduel de $S_i$.
\end{itemize}

Comme chaque $S_i$ est un quotient de $\CO_{\CX, \beta_i}^h$, la seconde propriété est bien satisfaite. Pour démontrer la première propriété, il suffit de montrer que le morphisme naturel $\CO_{\CX, \beta_{i-1}}^h \ra Frac(S_i)$ est surjectif.\\

Rappelons que par construction des points $\beta_i$, il existe un système régulier de paramètres $(f_1,...,f_i)$ de $\CO_{\CX,\beta_{i}}$ tel que $\beta_i = (f_1,...,f_i)$ et $\beta_{i-1} = (f_1,...,f_{i-1})$.
Posons $\tilde{\beta}_{i} \coloneqq \beta_{i} \CO_{\CX,\beta_i}^h$ et $\tilde{\beta}_{i-1} \coloneqq \beta_{i-1} \CO_{\CX,\beta_i}^h,$ l'idéal $\tilde{\beta}_{i-1}$ est bien premier puisque $S_i$ est un anneau intègre.\\
Appelons $\eta_i$ le point générique de $S_i$. Puisque l'inclusion $R_i \subseteq R_i^h = S_i$ est fidèlement plate (cf \cite{PMIHES_1967__32__5_0}, théorème 18.6.6, \textit{iii}), l'image de $\eta_i$ dans $Spec(\CO_{\CX,\beta_i})$ est le point $\beta_{i-1}$ par construction d'où le morphisme local $\CO_{\CX, \beta_{i-1}}^h \ra Frac(S_i)$. Considérons le diagramme commutatif suivant :

\[\begin{tikzcd}
 & \CO_{\CX,\beta_i}^h  & (\CO_{\CX,\beta_i}^h)_{\tilde{\beta}_{i-1}}\\
\CO_{\CX,\beta_i} &  & \text{Frac}(S_i)\\
\CO_{\CX,\beta_{i-1}} & \CO_{\CX,\beta_{i-1}}^h.\\
     \arrow[from=2-1, to=1-2, hook]
     \arrow[from=1-3, to=2-3, two heads]
     \arrow[from=1-2, to=1-3]
     \arrow[from=1-2, to=3-2, dotted]
     \arrow[from=2-1, to=3-1]
     \arrow[from=1-3, to=3-2]
     \arrow[from=3-2, to=2-3, dotted]
     \arrow[from=3-1, to=3-2]
\end{tikzcd}\]
dont les flèches pointillées proviennent de la propriété universelle de l'hensélisation. Précisons la construction du morphisme $(\CO_{\CX,\beta_i}^h)_{\tilde{\beta}_{i-1}} \ra \CO_{\CX,\beta_{i-1}}^h$ qui provient de la propriété universelle de la localisation. Soit donc $\alpha \in \CO_{\CX,\beta_i}^h \setminus \tilde{\beta}_{i-1}$. Si $\alpha \notin \tilde{\beta}_i$, alors $\alpha$ est inversible et il n'y a rien à faire. Supposons donc que $\alpha \in \tilde{\beta}_i \setminus \tilde{\beta}_{i-1}$. Par définition de $\beta_i$ et de $\beta_{i-1}$, on peut écrire $\alpha$ sous la forme $\alpha = \lambda f_i + \tilde{b}_{i-1}$ où $\lambda \in \CO_{\CX,\beta_i}^h$ et $\tilde{b}_{i-1} \in \tilde{\beta}_{i-1}$.\\

Vérifions dans un premier temps que l'image de $\tilde{\beta}_{i-1}$ dans $\CO_{\CX,\beta_{i-1}}^h$ est contenue dans l'idéal maximal $\beta_{i-1}\CO_{\CX,\beta_{i-1}}^h$. En effet, on a par définition que $\tilde{\beta}_{i-1} = \beta_{i-1}\CO_{\CX,\beta_i}^h$ et le diagramme précédent montre que tout élément de $\beta_{i-1}$ est bien envoyé dans $\beta_{i-1}\CO_{\CX,\beta_{i-1}}^h$. En particulier, l'image de $\tilde{b}_{i-1}$ n'est pas inversible.\\

Montrons ensuite que l'image de $\lambda f_i$ dans $\CO_{\CX,\beta_{i-1}}^h$ est inversible. D'une part, $f_i \in \beta_i \setminus \beta_{i-1}$ et on peut le relever en un élément de $\CO_{\CX,\beta_i}$ dont l'image dans $\CO_{\CX,\beta_{i-1}}^h$ est inversible par définition. D'autre part, on peut supposer, quitte à remplacer $f_i$ par une puissance de $f_i$, que $\lambda$ est inversible. En effet, si $\lambda \notin \tilde{\beta}_i$, $\lambda$ est déjà inversible et il n'y a rien à faire. Sinon, on a nécessairement que $\lambda \in \tilde{\beta}_i \setminus \tilde{\beta}_{i-1}$ et alors $\lambda = \lambda_1f_i + \tilde{b}_{1,n-1}$. En répétant ce processus tant que $\lambda_i$ n'est pas inversible, on obtient une suite $\{\lambda_i\}_{i\geq 0}$ telle que $\lambda_0 = \lambda$ et $\overline{\lambda_i}=f_i\overline{\lambda_{i+1}}$ dans $S_i$.
La suite des idéaux $<\overline{\lambda_i}>$  étant croissante, elle est donc nécessairement stationnaire : il existe donc un entier $k$ tel que $<\overline{\lambda_k}> = <\overline{\lambda_{k+1}}>$. Cela entraîne l'inversibilité de $\overline{f_i}$ ce qui est une contradiction puisque $\overline{f_i}$ est une uniformisante de l'anneau de valuation discrète $S_i$.\\

Enfin, comme l'image de $\lambda f_i$ est inversible mais que celle de $\tilde{b}_{i-1}$ ne l'est pas, on a que l'image de $\alpha$ est bien inversible car la somme de deux éléments non inversibles d'un anneau local n'est pas inversible.\\

Une chasse au diagramme montre que le morphisme $\CO_{\CX,\beta_{i-1}}^h \ra \text{Frac}(S_i)$ est surjectif, d'où l'isomorphisme voulu après avoir quotienté par $\beta_{i}$. \qed
\\
Penchons nous maintenant sur la cohomologie des points fermés.

\begin{lemme}\label{lemHarSzagénéral}
    Soit $k$ un corps de nombres et soit $Z$ une $k$-variété lisse géométriquement intègre de dimension $n$. Soit $\sF$ un faisceau étale sur $Z$. Alors, il existe un ensemble fini $S$ de places de $k$ contenant toutes les places archimédiennes tel que le morphisme 
$H^1(k,\sF) \ra \prod_{v \notin S}H^1(k_v,\sF)$
se factorise par le morphisme
$H^1(k,\sF) \ra \prod_{v \notin S}H^1(k_v,\sF)H^1(k,\sF) \ra \prod_{c \in Z^{(n)}}H^1(k(c),\sF).
$
En particulier, le groupe $$\Sha_n^1(k,\sF) \coloneqq Ker\bigg(H^1(k,\sF) \ra \prod_{c \in Z^{(n)}}H^1(k(c),\sF)\bigg)$$ s'injecte dans le groupe $$\Sha_{S}^1(k,\sF) \coloneqq Ker\bigg(H^1(k,\sF) \ra \prod_{v \notin S}H^1(k_v,\sF)\bigg).$$
\end{lemme}

\proof La preuve est identique à celle du lemme 2.1 de \cite{HarSza} qui traite le cas des courbes à l'aide des estimations de Lang-Weil \cite{LangWeil} et du lemme de Hensel.
\qed \\

On peut généraliser ce résultat au cas où la variété $Z$ n'est pas irréductible. Soit $Z = (\CX_k)$ la fibre spéciale de $\CX$ munie de sa structure réduite et soient $Z_1,...,Z_m$ les composantes irréductibles de $Z$. On pose $k_i \coloneqq k^s \cap k(Z_i)$.
\begin{lemme}\label{lemfinitude}
Le groupe $\Sha_{n}^1(k,G)$ est fini si et seulement si le noyau du morphisme
$$
H^1(k,G) \ra \prod_{i=1}^{m}H^1(k_i,G)
$$
est fini.
\end{lemme}

\proof Notons $n_i$ la dimension de la variété $Z_i$. Étant donné que chaque $Z_{i}$ est une $k_i$-variété géométriquement irréductible par normalité, le lemme \ref{lemHarSzagénéral} montre qu'il existe pour tout $i \in \{1,...,m\}$ un ensemble fini $S_i$ de places de $k_i$ contenant les places archimédiennes de $k_i$ tel que le noyau $\Sha^1_{n_i}(k_i,G)$ du morphisme  $$ H^1(k_i,G) \ra \prod_{c \in {Z_{i}}^{(n_i)}} H^1(k_i(c),G)$$ s'injecte dans le noyau $\Sha^1_{S_i}(k_i,G)$ du morphisme
$$
H^1(k_i,G) \ra \prod_{v \in \Omega_{k_i} \setminus S_i} H^1((k_i)_v,G)
$$
et ces derniers groupes sont de $m$-torsion finie pour tout $m>0$ d'après le lemme 2.2 de \cite{HarSza}. Par restriction-corestriction, les groupes $\Sha^1_{n_i}(k_i,G)$ sont donc également finis. En outre, on a le diagramme commutatif suivant :
\[\begin{tikzcd}
 H^1(k,G)  & \displaystyle{\prod_{i=1}^m} H^1(k_i,G)  & \\
 \displaystyle{\prod_{c \in Z^{(n)}}}H^1(k(c),G) & \displaystyle{\prod_{i=1}^m} \hspace{0,2cm}\displaystyle{\prod_{c \in {Z_{i}}^{(n_i)}}H^1(k_i(c),G)}
\end{tikzcd}\]
d'où l'injectivité du morphisme $$
\Sha^1_n(G)/ \text{Ker} \bigg(H^1(k, G) \ra \prod_{i=1}^m H^1(k_i, G)\bigg) \longrightarrow \displaystyle{\prod_{i=1}^m} \Sha^1_{n_i}(k_i,G).
$$
Il en résulte l'équivalence voulue. \qed \\

Rappelons que $X^{(1)}$ désigne le sous-ensemble de $\CX$ constitué des points de codimension $1$ correspondant à des points de $R$, i.e. ceux en dehors de la fibre singulière. Soit $\eta$ le point générique de $\CX$. Lorsque $R$ est de dimension 2, l'espace $X \coloneqq X^{(1)} \cup \{\eta\}$ est simplement l'ouvert $Spec( R) \setminus \{\fm\} $. Cependant, l'ensemble $X$ n'est plus muni d'une structure naturelle de schéma en dimension supérieure. Pour pallier ce problème, on définit pour tout entier $i \geq 0$ et tout faisceau étale $\sF$ défini sur un ouvert de $\CX$ contenant $X^{(1)}$ les groupes 
$$ 
H^i(X, \sF) \coloneqq \drl_{U \supseteq X^{(1)}} H^i(U,\sF)
$$
où la limite est prise sur les ouverts $U$ de $\CX$ contenant $X^{(1)}$. Il est préférable de travailler avec cet objet qui correspond moralement à la cohomologie d'un schéma de dimension 1.

Par abus de notation, on note également $G$ le groupe $G \times_k \CX$ . Pour établir la finitude du groupe $\Sha^1(K,G)$, on montre d'abord qu'il s'injecte dans l'image dans $H^1(K,G)$ du groupe
$$
\CK^1(X,G) \coloneqq Ker\bigg(H^1(X,G) \ra \prod_{c \in X^{(1)}} H^1(k(c), G)\bigg).
$$\\

Pour ce faire, on aura besoin d'une suite de localisation qui imitera celle existante pour les schémas de Dedekind, cf le lemme 2.3 de \cite{HarSza}.

\lemme\label{lemlocalisation} Soit $S=\{s_1,...,s_n\}$ un sous-ensemble fini de $X^{(1)}$. Alors, le noyau du morphisme de restriction
$$
\drl_{U \supseteq X^{(1)}} H^i(U\setminus (\ol{S} \cap U),G) \ra \bigoplus_{v \in S} H^i(K_v,G)
$$

est inclus dans l'image du morphisme de restriction
$$
\drl_{U \supseteq X^{(1)}} H^i(U,G) \ra \drl_{U \supseteq X^{(1)}} H^i(U\setminus (\ol{S} \cap U),G).
$$

\proof Soit $U \supseteq X^{(1)}$ un ouvert de $\CX$, la suite de localisation pour la cohomologie étale fournit la suite exacte
$$
H^i(U,G) \ra H^i(U\setminus (\ol{S} \cap U),G) \ra H^i_{\ol{S}\cap U}(U,G)
$$
qui, après passage à la limite directe, donne
$$
\drl_{U \supseteq X^{(1)}} H^i(U,G) \ra \drl_{U \supseteq X^{(1)}} H^i(U\setminus (\ol{S} \cap U),G) \ra \drl_{U \supseteq X^{(1)}} H^i_{\ol{S}\cap U}(U,G).
$$
Remarquons d'abord que l'on peut se ramener au cas où $S=\{s\}$ est un singleton. En effet, quitte à ne considérer que des ouverts $U$ suffisamment petits, on peut supposer que les ensembles $\ol{\{s_i\}} \cap U$ ne s'intersectent pas et on alors par la suite de Mayer-Vietoris pour la cohomologie à support que 
$$
H^i_{\ol{S}\cap U}(U,G) \simeq \bigoplus_{s_i \in S} H^i_{\ol{\{s_i\}}\cap U}(U,G).
$$
\\
Construisons en premier lieu un isomorphisme entre les groupes $ \displaystyle \drl_{U \supseteq X^{(1)}}H^i_{\ol{\{s\}}\cap U}(U,G)$ et $H^i_s(\CO_s^h,G)$.
Soit donc un ouvert $U \supseteq X^{(1)}$. D'une part, en passant à la limite sur tous les voisinages étales élémentaires $(Y,y) \ra (U,s)$ avec $Y$ affine tel que $f^{-1}(s)=\{y\}$, on a par le corollaire 5.8, exposé $\textit{VII},$ de \cite{SGA4} un isomorphisme
$$
\drl_{(Y,y)}H^i(Y,G) \simeq H^i(\CO_s^h,G).
$$
 Considérons le diagramme commutatif suivant provenant des suites de localisation :

\[\begin{tikzcd}
 ...&\drl_{(Y,y)} H^i(Y,G)  &\drl_{(Y,y)} H^i(Y \setminus \ol{\{y\}}, G) &\drl_{(Y,y)} H^i_{\ol{\{y\}}}(Y,G) & ....\\
... & H^i(\CO_s^h,G) & H^i(K_s^h,G) & H^i_s(\CO_s^h,G) & ...\\
     \arrow[from=1-1, to=1-2]
     \arrow[from=1-2, to=1-3]
     \arrow[from=1-3, to=1-4]
     \arrow[from=1-4, to=1-5]
     \arrow[from=2-1, to=2-2]
     \arrow[from=2-2, to=2-3]
     \arrow[from=2-3, to=2-4]
     \arrow[from=2-4, to=2-5]
     \arrow[from=1-2, to=2-2, "\simeq"]
     \arrow[from=1-3, to=2-3, "\simeq"]
     \arrow[from=1-4, to=2-4]
\end{tikzcd}\]

Le lemme des cinq montre que l'on dispose d'un isomorphisme
$$
\drl_{(Y,y)}H^i_{\ol{\{y\}}}(Y,G) \simeq H^i_s(\CO_s^h,G).
$$

D'autre part, si l'on considère un tel voisinage étale $(Y,y) \ra (U,s)$, il induit un morphisme birationnel $\ol{\{y\}} \ra \ol{\{s\}} \cap U$ qui, quitte à réduire $U$, est un isomorphisme. Il en résulte par le lemme d' excision un isomorphisme
$$
\drl_{(Y,y)}H^i_{\ol{y}}(Y,G) \simeq \drl_{U \supseteq X^{(1)}} H^i_{\ol{S}\cap U}(U,G).
$$
On obtient au final que
$$
\drl_{U \supseteq X^{(1)}} H^i_{\ol{S}\cap U}(U,G) \simeq H^i_s(\CO_s^h,G)
$$
et on peut alors conclure en utilisant la suite de localisation sur $Spec(O_s^h)$. On réécrit ici l'argument provenant de la preuve du lemme 2.3 de \cite{HarSza}) : on considère la suite
$$
H^i(\CO_s^h,G) \ra H^i(K_s^h,G) \ra H^{i+1}_s(\CO_s^h,G).
$$
On dispose d'un isomorphisme $H^i(K_s^h,G) \simeq H^i(K_s,G)$ par le théorème d'approximation de Greenberg \cite{Greenberg1966} (cf le lemme 2.7 de \cite{HarSza3} pour une preuve de cet isomorphisme). Considérons le diagramme commutatif suivant :

\[\begin{tikzcd}
     &\displaystyle \drl_{U \supseteq X^{(1)}} H^i_{\ol{S}\cap U}(U,G) & H^i_s(\CO_s^h,G) \\
 H^i(\CO_s^h,G) &  H^i(K_s,G) & H^{i+1}_s(\CO_s^h,G).
     \arrow[from=1-2, to=1-3]
     \arrow[from=1-3, to=2-3, "\simeq"]
     \arrow[from=2-1, to=2-2]
     \arrow[from=2-2, to=2-3]
     \arrow[from=1-2, to=2-2]
\end{tikzcd}\]

On en déduit le résultat voulu.

\qed

L'objet du lemme suivant est de montrer que ce résultat nous permet de ramener l'étude du groupe $\Sha^1(K,G)$ à celle du groupe $\CK^1(X,G)$.

\lemme\label{lemfinal} Le groupe $\Sha^1(K,G)$ est contenu dans l'image de $\CK^1(X,G)$ dans $H^1(K,G)$.

\proof Soit $\alpha \in \Sha^1(K,G)$. On peut relever $\alpha$ en un élément de $H^1(V,G)$ pour un ouvert $V$ de $\CX$. Quitte à réduire $V$, on peut supposer que les seuls points de codimension $1$ de $\CX$ appartenant à $V$ sont des éléments de $X^{(1)}$. Soient alors $s_1,...,s_i,t_1,...,t_j$ les points génériques des composantes irréductibles de $\CX \setminus V$ où $s_1,...,s_i$  sont des éléments de $X^{(1)}$ tandis que $t_1,...,t_j \notin X^{(1)}$ 
L'ensemble $U \coloneqq (\CX \setminus \ol{\{t_1,...,t_j\}} )$ est un ouvert de $\CX$ tel que $\CX \setminus V = \ol{\{s_1,...,s_i,t_1,...,t_j\}} = \ol{\{s_1,...,s_i\}} \cup \ol{\{t_1,...,t_j\}} = (\CX \setminus U) \cup \ol{S}$ par définition.

Ainsi, on a bien que $V = U \setminus (\ol{S} \cap U)$ et $U \supseteq X^{(1)}$; le lemme \ref{lemlocalisation} montre donc que $\alpha$ se relève à fortiori en un élément de $H^1(X,G)$. Soit maintenant $v \in X^{(1)}$, on a par lisseté de $G$ que $H^1(O_v^h,G) \simeq H^1(k(v),G)$ (cf \cite{Milne}, \textit{III}, remarque 3,11,a). De plus, puisque l'anneau $O_v^h$ est un anneau de valuation discrète, le morphisme naturel $H^1(O_v^h,G) \ra H^1(K_v^h,G)$ est injectif (cf le théorème 4.5 de \cite{Nisnevich}). On conclut en notant que $H^1(K_v^h,G) \simeq H^1(K_v,G)$ par le théorème d'approximation de Greenberg \cite{Greenberg1966}.
\qed

\begin{prop}\label{prop1}
Les groupes 
$$
\Sha^1_1(k,G) \coloneqq Ker\bigg(H^1(k,G) \ra \prod_{{c \in X^{(1)}}}H^1(k(c),G) \bigg)
$$
et 
$$
\Sha^1_n(k,G) \coloneqq Ker\bigg(H^1(k,G) \ra \prod_{c \in \CX^{(n)}}H^1(k(c),G) \bigg)
$$
sont égaux.
\end{prop}

Insistons sur le fait que l'on ne considère pas tous les points de codimension $1$ dans la définition de $\Sha^1_1(k,G)$.

\proof On se contente de montrer l'égalité $\Sha^1_1(k,G) = \Sha^1_n(k,G)$, le cas général étant identique. Remarquons tout d'abord que pour tout  $\beta_1 \in X^{(1)}$, il existe un point $\beta_n \in \CX^{(n)}$ tel que $\beta_n \in \overline{\{\beta_1\}}$ par le lemme \ref{lemdimension}. De plus, on a par (\cite{Milne}, \textit{III}, remarque 3,11,a) des isomorphismes $H^1(k(\beta_n),G) \simeq H^1(\CO_{\CX,\beta_n}^h,G)$ et $H^1(\CO_{\CX,\beta_1}^h,G) \simeq H^1(k(\beta_1),G).$

L'inclusion $\CO_{\CX,\beta_n} \ra \CO_{\CX,\beta_1}$ induit un morphisme $\CO_{\CX,\beta_n}^h \ra \CO_{\CX,\beta_1}^h$ qui fournit en prenant la cohomologie un morphisme
$$
H^1(k(\beta_n),G) \simeq H^1(\CO_{\CX,\beta_n}^h,G) \ra H^1(\CO_{\CX,\beta_1}^h,G) \simeq H^1(k(\beta_1),G).
$$
Cela montre que le morphisme $H^1(k,G) \ra H^1(k(\beta_1),G)$ se factorise par $H^1(k(\beta_n),G) \ra H^1(k(\beta_1),G)$. On en déduit une première inclusion  $\Sha^1_n(k,G) \subseteq \Sha^1_1(k,G)$.\\

Soit maintenant $\alpha \in \Sha^1_1(k,G)$ et soit $\beta_n$ un point de codimension $n$ de $\CX$. D'après le lemme \ref{lemdrapeau} il existe un modèle régulier $\CX'$ de $\CX$ et un drapeau régulier $(\beta_1,...,\beta_n)$ dans $\CX'$ tel que $\beta_1 \in X^{(1)}$ et de sorte que pour tout $i \in \{2,...,n\}$, il existe des anneaux $S_i$ henséliens et réguliers satisfaisant $Frac(S_{i+1}) \simeq k(\beta_i)$ et $k(S_i) \simeq  k(\beta_i)$.
Il ne reste plus qu'à réutiliser l'argument du lemme \ref{lemfinal} : le morphisme naturel $H^1(S_{i+1},G) \ra H^1(k(\beta_i),G)$ est injectif par régularité de $S_{i+1}$ et $H^1(S_{i+1},G) \simeq H^1(k(\beta_{i+1}),G)$. Itérer ce processus $n-1$ fois montre que l'image de $\alpha$ dans $H^1(k(\beta_n),G)$ est nulle, c'est ce que l'on voulait démontrer.
\qed \\

\rmke On peut en fait montrer que pour tout entier $1 \leq i \leq n$, les groupes
$$
\Sha^1_i(k,G) \coloneqq Ker\bigg(H^1(k,G) \ra \prod_{\underset{\text{codim}(f(c)) = i}{c \in \CX^{(i)},}}H^1(k(c),G) \bigg)
$$
sont égaux en utilisant le lemme \ref{lemdrapeau}.\\

Terminons cette section avec un petit lemme qui servira plusieurs fois par la suite.

\begin{lemme}\label{lemzariski}
    Soit $X \ra S$ un morphisme propre birationnel de schémas noethériens avec $X$ intègre et $S$ normal. Soit $U$ un ouvert de $S$ contenant tous les points de codimension $1$. Alors, $H^0(X,\CO_X) \simeq H^0(S,\CO_S) \simeq H^0(U,\CO_U)$.
\end{lemme}

\proof 
D'après le théorème principal de Zariski (cf \cite{Hartshorne},\textit{III}, corollaire 11.4), on dispose d'un isomorphisme de faisceaux $f_*\CO_X = \CO_S$ qui donne le premier isomorphisme après l'évaluation en $S$.
Le second isomorphisme provient du lemme de Hartog algébrique, cf le théorème 1.14 du chapitre $4$ de \cite{Liu}. \qed

\section{Le cas déployé}\label{secdéployé}
Afin d'étudier le groupe $\CK^1(X,G)$, on traite d'abord le cas où $G$ est déployé afin de pouvoir utiliser les propriétés des groupes de Picard. Pour ce faire, soit $l/ k$ une extension finie déployant $G$ et notons $R_l \coloneqq R \times_k l$, $\CY \coloneqq \CX \times_R R_l$. L'anneau $R_l$ est un $R$-module fini, il est donc également noethérien, local de corps résiduel $l$, complet, excellent et intègre par hypothèse sur $R$.  Par le théorème d'Hilbert 90, le groupe $\CK^1(Y,G)$ est isomorphe à un produit fini de $H^1(Y, \G_m)$ et de $\CK^1(Y, \Z / n\Z)$, on se contentera donc d'étudier des groupes de cette forme.

\subsection{Les tores}\label{sectiontores}
Le but de cette section est de montrer que le groupe $\ol{\Pic(\CY)}$  est de type fini, où $\ol{\Pic(\CY)}$ désigne le quotient de $\Pic(\CY)$ par son sous-groupe divisible maximal.
Notons $\fm_l$ l'idéal maximal de $R_l$. Pour tout $n \geq 0$, notons $Y_n$ la fibre de l'idéal $\fm_l^{n+1}$.

D'après la scholie 5.1.7 de \cite{PMIHES_1961__11__5_0}, on dispose d'un isomorphisme 
$$
\Pic(\CY) \simeq \prl \Pic(Y_n).
$$ Cela nous incite à essayer de comprendre le groupe $\Pic(\CY)$ à travers les groupes $\Pic(Y_n)$.\\

À cet effet, on se penche d'abord sur le groupe de Picard de la fibre spéciale munie de sa structure réduite $Z \coloneqq (Y_0)_{red}$.

\begin{lemme}\label{lemPicfinitude}
 Soit $Z = (Y_0)_{red}$ la fibre spéciale de $\CY$ munie de sa structure réduite et soient $Z_1,...Z_n$ les composantes irréductibles de $Z$. Supposons que l'hypothèse \ref{hyp2} soit satisfaite. Alors, le groupe $\Pic(Z)$ est de type fini.
\end{lemme}

\proof
La fibre spéciale de $\CY$ étant un diviseur à croisements normaux strict, $Z$ est une  $l$-variété propre dont les composantes irréductibles $Z_1,...,Z_l$  sont lisses. D'après le théorème de Mordell-Weil et le théorème de Néron-Severi, les groupes $\Pic^0(Z_i)$ et $NS(Z_i)$ sont de type fini; il en va donc de même pour chaque $\Pic(Z_i)$ d'après la suite exacte
$$
0 \ra \Pic^0(Z_i) \ra \Pic(Z_i) \ra NS(Z_i) \ra 0.
$$

Par l'hypothèse \ref{hyp2}, le noyau du morphisme $Pic(Z) \ra \prod_{i=1}^n \Pic(Z_i)$ est fini; on en conclut que le groupe $\Pic(Z)$ est de type fini. \qed 
\\

\rmke Intéressons-nous aux situations où le noyau du morphisme $\Pic(Z) \ra \prod_{i=1}^n \Pic(Z_i)$ est fini. Lorsque $R$ est de dimension $2$, la fibre spéciale est une union de courbes régulières et on peut considérer la suite exacte suivante dont l'exactitude se déduit localement :
$$
0 \ra \CO_Z^* \ra \prod_{i=1}^n \CO_{Z_i}^* \ra \bigoplus_{p \in P} (i_p)_*\gm \ra 0
$$
où $P$ désigne l'ensemble des points d'intersections des composantes irréductibles de $Z$ et $(i_p)_*$ désigne le faisceau gratte-ciel concentré en le point $p \in P$. Cette suite induit la suite exacte cohomologique suivante
$$
\prod_{i=1}^n H^0(Z_i, \gm) \ra \prod_{p \in P} H^0(k(p), \gm)  \ra \Pic(Z) \ra \prod_{i=1}^n \Pic(Z_i).
$$
Ainsi, il suffit que le conoyau du morphisme $\prod_{i=1}^n H^0(\CO_{Z_i}, \gm) \ra \prod_{p \in P} H^0(k(p), \gm)$ soit fini. C'est notamment le cas lorsque le graphe associé à $Z$ est un arbre et que les points d'intersection des composantes irréductibles de $Z$ sont rationnels; cette situation apparaît par exemple quand la singularité de $Spec(R)$ est rationnelle, cf \cite{Lipman-rational}.
De manière plus générale, si la fibre spéciale de $\CX$ est irréductible, alors elle est également géométriquement irréductible puisque elle est géométriquement connexe par le lemme \ref{lemstein}. \\

On supposera donc pour le reste de l'article que l'hypothèse \ref{hyp2} est satisfaite. La prochaine étape du raisonnement est de comparer les groupes $\Pic(Y_0)$ et $\Pic({(Y_0)}_{red}).$

\begin{lemme}\label{Picred}
Le morphisme naturel
    $$\Pic(Y_0) \ra \Pic({(Y_0)}_{red}) $$ induit un isomorphisme entre $\ol{\Pic(Y_0)}$ et un sous-groupe de $\Pic({{(Y_0)}_{red}})$. En particulier, $\ol{\Pic(Y_0)}$ est de type fini.\\
    De plus, le sous-groupe divisible maximal de $\Pic(Y_0)$ est sans torsion.
\end{lemme}

\proof Notons $r : {(Y_0)}_{red} \ra Y_0$ l'immersion fermée naturelle et $\CI$ le faisceau d'idéal quasi-cohérent pour la topologie étale correspondant au noyau de $r$.
Considérons la suite exacte de faisceaux étales sur $Y_0$ :
$$  0 \ra \CI \ra \G_{m,{Y_0}} \ra r_*\G_{m,{(Y_0)}_{red}} \ra 0 .$$

Prouvons en premier lieu que le morphisme induit $H^0(Y_0, \G_{m}) \ra H^0({(Y_0)}_{red}, \G_{m}) $ est surjectif : en effet, le schéma ${(Y_0)}_{red}$ est une $l-$variété propre, géométriquement connexe par le lemme \ref{lemstein} et (géométriquement) réduite ce qui implique que $H^0({(Y_0)}_{red}, \G_{m}) = l^*$. Ainsi, on dispose de morphismes 
$$
l^* \ra H^0(Y_0, \G_{m}) \ra H^0({(Y_0)}_{red}, \G_{m}) = l^*
$$
dont la composition est l'identité; il en résulte la surjectivité voulue.

De ces considérations résultent la suite exacte de cohomologie 
$$
 0 \ra  H^1(Y_0,\CI) \ra \Pic(Y_0) \ra H^1(Y_0, {r}_*(\gm)_{{(Y_0)}_{red}}) \ra H^2(Y_0,\CI) .
$$

Le groupe $H^1(Y_0,\CI)_{\textit{\'et}} \simeq H^1(Y_0,\CI)_{Zar}$ est uniquement divisible d'où un isomorphisme $$\ol{\Pic(Y_0)} \simeq  \ol{Ker\big({H^1(Y_0, {r}_*\G_{m})} \ra H^2(Y_0,\CI) \big)}.$$ \\
Considérons maintenant la suite spectrale de Leray
$$
H^i(Y_0, R^j{r}_*(\G_m)) \Ra H^{i+j}({(Y_0)}_{red}, \G_m)
$$
Notons d'abord que $R^1{r}_*(\gm)=0$ : en effet, pour tout point géométrique $\ol{x}$, on a un isomorphisme $R^1{r}_*(\gm)_{\ol{x}} = H^1({(Y_0)}_{red} \times Spec(\CO_{Y_O,\ol{x}}, \gm))$ et $r$ étant une immersion fermée, le schéma ${(Y_0)}_{red} \times_{\ol{\{x\}}} Spec(\CO_{Y_0,\ol{x}})$ est isomorphe au spectre d'un quotient de $\CO_{Y_0,\ol{x}}$ qui est également un anneau local. Sachant que le groupe de Picard d'un anneau local est trivial (\cite{Milne},\textit{III}, lemme 4.10), on en déduit la nullité de $R^1{r}_*(\G_m)$.\\
La suite exacte à cinq termes fournit alors un isomorphisme
$$H^1(Y_0, r_*(\G_m)) \simeq H^1({(Y_0)}_{red}, \G_m) = \Pic({(Y_0)}_{red}).$$

Le groupe $\Pic({(Y_0)}_{red})$ est de type fini par le lemme \ref{lemPicfinitude} et n'a donc aucun sous-groupe divisible non trivial, il en résulte que
$$
\ol{\Pic(Y_0)} \simeq Ker\big({H^1(Y_0, {r}_*\G_{m})} \ra H^2(Y_0,\CI) \big) \hookrightarrow \Pic((Y_0)_{red}),
$$
c'est ce que l'on voulait démontrer.
\qed
\\

\rmke Lorsque l'anneau $R$ est de dimension 2, la fibre spéciale est de dimension 1. Comme le faisceau $\CI$ est cohérent, on peut passer à la cohomologie de Zariski ce qui implique la nullité de $H^2(Y_0,\CI)$ d'où un isomorphisme $\ol{\Pic(Y_0)} \simeq \Pic({(Y_0)}_{red})$.\\

On souhaite maintenant relier le groupe $\ol{\Pic(\CY)}$ au groupe $\ol{\Pic(Y_0)}$.

\begin{lemme}\label{PicYn} Rappelons que $Y_n$ désigne la fibre de $\fm^{n+1}$ pour tout entier $n$.
Alors, $\ol{\Pic(Y_{n+1})} \hookrightarrow \ol{\Pic(Y_n)}$ et
$$\prl_{n \geq 0} \ol{\Pic(Y_n)} = \bigcap_{n \geq 0} \ol{\Pic(Y_n)}.$$
De plus, le sous-groupe divisible maximal de $\Pic(\CY)$ est sans torsion.
\end{lemme}

\proof  Considérons la suite exacte de faisceaux de Zariski sur $\CY$ 

$$ 0 \ra \frac{\fm^n_l \CO_\CY }{\fm^{n+1}_l \CO_Y} \ra \bigg(\frac{\CO_\CY }{\fm^{n+1}_l \CO_\CY}\bigg)^* \ra \bigg(\frac{\CO_\CY }{\fm^{n}_l \CO_\CY}\bigg)^* \ra 1.$$
Passer à la cohomologie donne la suite exacte suivante :
\begin{equation}\label{suiteexactePic}
 H^0(Y_{n-1}, \gm) \ra H^1\bigg(\CY,\frac{\fm^n_l \CO_\CY }{\fm^{n+1}_l}\bigg) \ra \Pic(Y_n) \ra \Pic(Y_{n-1}) \ra H^2\bigg(\CY,\frac{\fm^n_l \CO_\CY }{\fm^{n+1}_l}\bigg).
\end{equation}

D'après le lemme \ref{lemzariski}, on a que $\CO_{\CY}(\CY)= \CO_{R_l}(R_l) = R_l$.
L'anneau $R_l/\fm^{n+1}_l$ étant hensélien, le morphisme $H^0\bigg(\CY,\bigg(\frac{\CO_\CY }{\fm^{n+1}_l \CO_\CY}\bigg)^*\bigg)  \ra H^0\bigg(\CY,\bigg(\frac{\CO_\CY }{\fm^{n}_l \CO_\CY}\bigg)^*\bigg)$ est surjectif ce qui implique par la suite exacte longue \ref{suiteexactePic} l'injectivité du morphisme $H^1(\CY,\frac{\fm^n_l \CO_\CY }{\fm^{n+1}_l}) \ra \Pic(Y_n).$ \\
 Il en résulte la suite exacte suivante
 $$
 0 \ra H^1(\CY,\frac{\fm^n_l \CO_\CY }{\fm^{n+1}_l}) \ra \Pic(Y_n) \ra Ker\bigg( \Pic(Y_{n-1}) \ra H^2(\CY,\frac{\fm^n_l \CO_\CY }{\fm^{n+1}_l}) \bigg) \ra 0.
 $$
Supposons que $n=1$. D'après le lemme \ref{Picred}, $\Pic(Y_0) \simeq U_0 \bigoplus F_0$ où $U_0$ est un groupe uniquement divisible et $F_0$ est un groupe de type fini (qui est un sous groupe de $\Pic({(Y_0)}_{red})$. 
 Le groupe $H^2(\CY,\frac{\fm^n_l \CO_Y }{\fm^{n+1}_l})$ est également un groupe uniquement divisible et le noyau d'un morphisme entre deux groupes uniquement divisibles est à son tour un groupe uniquement divisible, aussi le noyau du morphisme $\varphi : \Pic(Y_0) \ra H^2(\CY,\frac{\fm_l \CO_Y }{\fm^{2}_l})$ s'écrit comme une somme directe entre un groupe uniquement divisible $U'_0 \subset U_0$ et un groupe de type fini $F'_0 \subset F_0$.\\
En outre, le groupe $H^1(\CY,\frac{\fm_l \CO_Y }{\fm^{2}_l})$ est un $l$-espace vectoriel, c'est donc un groupe uniquement divisible. Il est à fortiori injectif ce qui nous permet d'obtenir la suite exacte scindée suivante 
 $$
0 \ra H^1\bigg(\CY,\frac{\fm_l \CO_\CY }{\fm^{2}_l}\bigg) \ra \Pic(Y_1) \ra U'_0 \oplus F'_0 \ra 0.
 $$

Ainsi, $\Pic(Y_1)\simeq U_1 \bigoplus F'_0$ où $U_1 \coloneqq H^1\bigg(\CY,\frac{\fm_l \CO_\CY }{\fm^{2}_l}\bigg) \bigoplus U'_0$ est un groupe uniquement divisible; cela implique en plus que $\ol{\Pic(Y_1)} \hookrightarrow \ol{\Pic(Y_0)}$. On peut alors itérer le processus et en conclure que $\ol{\Pic(Y_{n+1})} \hookrightarrow \ol{\Pic(Y_n)}$ pour tout $n$ ce qui fournit immédiatement l'égalité

$$\prl_{n \geq 0} \ol{\Pic(Y_n)} = \bigcap_{n \geq 0} \ol{\Pic(Y_n)}.$$
 
 \qed

\rmke On obtient ici aussi un résultat plus précis en dimension 2. Le faisceau $\frac{\fm^n_l \CO_Y }{\fm^{n+1}_l \CO_Y}$ étant quasi-cohérent, on peut passer à la cohomologie de Zariski où le théorème de changement de base propre combiné au fait que la fibre $Y_0$ soit de dimension un montrent que le groupe $H^2(\CY,\frac{\fm^n_l \CO_Y }{\fm^{n+1}_l}) = H^2(X_0,\frac{\fm^n_l \CO_Y }{\fm^{n+1}}) $ est trivial.
Le même raisonnement que celui employé dans la preuve ci-dessus fournit un isomorphisme entre $\ol{\Pic(Y_n)}$ et $\ol{\Pic(Y_{n+1})}$ pour tout $n \geq 0$, d'où un isomorphisme $\prl_{n \geq 0} \ol{\Pic(Y_n)} \simeq \ol{\Pic(Y_0)}$.\\

L'objet du prochain lemme est de relier la structure des groupes de Picard des épaississements de la fibre spéciale de $\CY$ avec celle du groupe de Picard de $\CY$.

\begin{lemme}\label{lemPicfini}
Le groupe $\ol{\Pic(\CY)}$ est de type fini. De plus, le sous-groupe divisible maximal de $\Pic(\CY)$ est sans torsion.
\end{lemme}

\proof
Considérons les suites exactes
$$
0 \ra D_n \ra \Pic(Y_n) \ra \ol{\Pic(Y_n)}  \ra 0
$$
où $D_n$ désigne le sous-groupe divisible maximal de $\Pic(Y_n)$. D'après le lemme \ref{PicYn}, ce groupe est en fait uniquement divisible. Comme l'image d'un groupe divisible est divisible, ces suites exactes forment un système projectif. Par conséquent,après être passé à la limite, on obtient la suite exacte suivante :
$$
0 \ra \prl  D_n \ra  \prl \Pic(Y_n) \ra \prl \ol{\Pic(Y_n)}. 
$$

Par la scholie 5.1.7 de \cite{PMIHES_1961__11__5_0}, on dispose d'un isomorphisme 
$$
\Pic(\CY) \simeq \prl \Pic(Y_n).
$$

Enfin, puisque le groupe $\prl  D_n$ est uniquement divisible (en tant que limite de groupes uniquement divisibles) et que $\bigcap_{n \geq 0} \ol{\Pic(Y_n)}$ est de type fini (et donc sans sous-groupe divisible non trivial), $\prl  D_n$ est le sous-groupe divisible maximal de $\Pic(\CY)$ d'où l'injectivité du morphisme $\ol{\Pic(\CY)} \ra \bigcap_{n \geq 0} \ol{\Pic(Y_n)}.$ On déduit finalement que le groupe $\ol{\Pic(\CY)}$ est bien de type fini. \qed\\

La dernière étape consiste à étudier la structure du groupe
$$\Pic(Y) \coloneqq \drl_{U \supseteq Y^{(1)}} \Pic(U)$$
à l'aide des résultats précédemment obtenus.\\

\begin{prop}\label{propPicfini}
    Le groupe 
$\ol{\Pic(Y)}$ est de type fini. De plus, la $m$-torsion du sous-groupe divisible maximal de $Pic(Y)$ est finie pour tout entier $m>0$. 
\end{prop}
\proof Pour tout ouvert $U \supseteq \CY^{(1)}$ dont le complémentaire est de codimension $1$, notons $n_U$ le nombre de diviseurs de $\CY$ à support dans $\CY \setminus U$. Par régularité de $\CY$, on dispose d'une suite exacte
$$
\bigoplus_{i=1}^{n_U} \Z \ra \Pic(\CY) \ra \Pic(U) \ra 0.
$$

Quitte à ne considérer que les ouverts $U \supseteq Y^{(1)}$ inclus dans le locus régulier de $Spec(R)$, on peut supposer que $n_U=n$ pour tout $U$, l'entier $n$ désignant le nombre de diviseurs à support dans la fibre singulière.
Les morphismes de restriction induisant des suites exactes compatibles, on peut passer à la limite directe et obtenir la suite exacte
$$
\bigoplus_{i=1}^n \Z \ra \Pic(\CY) \ra \Pic(Y) \ra 0.
$$
Il en résulte que le groupe $\ol{\Pic(Y)}$ est bien de type fini.\\

Montrons maintenant que cette dernière suite exacte prouve également que la $m$-torsion du sous groupe divisible maximal de $\Pic(Y)$ est finie pour tout entier naturel $m>0$. 
Le lemme du serpent appliqué au diagramme
\[\begin{tikzcd}
   \bigoplus_{i=1}^n \Z  & \Pic(\CY) & \Pic(Y) & 0 \\
   \bigoplus_{i=1}^n \Z  & \Pic(\CY) & \Pic(Y) & 0\\
    \arrow[from=1-1, to=2-1, "\times m"]
	\arrow[from=1-1, to=1-2]
    \arrow[from=1-2, to=1-3]
    \arrow[from=1-3, to=1-4]
	\arrow[from=2-1, to=2-2]
    \arrow[from=2-2, to=2-3]
    \arrow[from=2-3, to=2-4]
    \arrow[from=1-2, to=2-2, "\times m"]
    \arrow[from=1-3, to=2-3, "\times m"]
    \arrow[from=1-4, to=2-4]
\end{tikzcd}\]
donne la suite exacte
$$
 \Pic(\CY)_m \ra \Pic(Y)_m \ra A/mA
$$
où $A$ est un groupe de type fini. Comme le groupe $\ol{\Pic(\CY)}$ est également de type fini et que son sous-groupe divisible maximal est sans torsion par le lemme \ref{lemPicfini}, les groupes ${}_{m}\Pic(\CY)$ et $A/mA$ sont tous deux finis et il en va donc de même du groupe ${}_{m}\Pic(Y)$. On en déduit le résultat voulu. \qed

\subsection{Les groupes finis}
Établissons maintenant un résultat de finitude analogue lorsque les coefficients sont des groupes abéliens finis.

\begin{prop}\label{propdéployéfinitude}
    Soit $m > 0$. Alors, le groupe $\CK^1(\CY, \Z / m\Z)$ est fini.
\end{prop}

\proof
Notons premièrement que pour tout ouvert $U \supseteq Y^{(1)}$, le groupe $H^1(U, \Z / m\Z)$ s'inscrit dans une suite exacte
$$
0 \ra H^1(l,\Z / m\Z) \ra H^1(U, \Z / m\Z) \ra H^1(U_s, \Z / m\Z)
$$
où $U_s$ désigne le changement de base de $U$ à une clôture algébrique $l_s$ de $l$. Passer à la limite directe donne la suite exacte
$$
0 \ra H^1(l,\Z / m\Z) \ra \drl_{U \supseteq Y^{(1)}} H^1(U, \Z / m\Z) \ra \drl_{U \supseteq Y^{(1)}} H^1(U_s, \Z / m\Z).
$$
Établissons en premier lieu la finitude du groupe $\drl_{U \supseteq Y^{(1)}} H^1(U_s, \Z / m\Z)$. Remarquons d'abord que l'on peut se contenter de considérer les ouverts $U \supseteq Y^{(1)}$ inclus dans le lieu régulier de $Spec(R_l)$. Soit donc $U$ un tel ouvert; la suite de Kummer fournit le diagramme commutatif
\[\begin{tikzcd}
   0 & \CO_{\CY_s}(\CY_s)^*/m & H^1(\CY_s, \mu_m) & {}_{m}\Pic(\CY_s) & 0 \\
	0 & \CO_{U_s}(U_s)^*/m & H^1(U_s, \mu_m) & {}_{m}\Pic(U_s) & 0.\\
	\arrow[from=1-1, to=1-2]
    \arrow[from=1-2, to=1-3]
    \arrow[from=1-3, to=1-4]
    \arrow[from=1-4, to=1-5]
	\arrow[from=2-1, to=2-2]
    \arrow[from=2-2, to=2-3]
    \arrow[from=2-3, to=2-4]
    \arrow[from=2-4, to=2-5]
    \arrow[from=1-2, to=2-2]
    \arrow[from=1-3, to=2-3]
    \arrow[from=1-4, to=2-4]
\end{tikzcd}\]

La première flèche verticale est un isomorphisme par la proposition \ref{prop1}. La troisième flèche verticale a pour conoyau un groupe fini de cardinal borné par $m^n$ où $n$ désigne le nombre de composantes irréductibles de la fibre singulière $\CX_{sing}$ : en effet, cela provient du lemme du serpent appliqué au diagramme
\[\begin{tikzcd}
    \Z^{n} & \Pic(\CY_s) & \Pic(U_s) & 0 \\
    \Z^{n} & \Pic(\CY_s) & \Pic(U_s) & 0\\
    \arrow[from=1-1, to=2-1, "\times m"]
	\arrow[from=1-1, to=1-2]
    \arrow[from=1-2, to=1-3]
    \arrow[from=1-3, to=1-4]
	\arrow[from=2-1, to=2-2]
    \arrow[from=2-2, to=2-3]
    \arrow[from=2-3, to=2-4]
    \arrow[from=1-2, to=2-2, "\times m"]
    \arrow[from=1-3, to=2-3, "\times m"]
    \arrow[from=1-4, to=2-4]
\end{tikzcd}\]
par régularité de $\CY_s$.\\

Ainsi, le conoyau du morphisme $H^1(\CY_s, \mu_m) \ra H^1(U_s, \mu_m)$ est lui aussi fini par le lemme du serpent et son cardinal est borné par $m^n$ pour tout ouvert $U_s$. D'après le corollaire 5.2 de \cite{SGA4}, \textit{XIX}, la groupe $H^1(\CY_s, \mu_m)$ est fini; cela implique que les groupes $H^1(U_s, \mu_m)$ sont également finis pour tout ouvert $U_s$ contenant $Y^{(1)}$ et le cardinal de tous les groupes $H^1(U_s, \mu_m)$ est borné par $m^n*card(H^1(\CY_s, \mu_m))$.\\

Montrons à présent que cela implique la finitude du groupe $\drl_{U \supseteq Y^{(1)}} H^1(U_s, \Z / m\Z)$ : en effet,  supposons par l'absurde que ce groupe soit infini et soient $\alpha_1,...,\alpha_M$ des éléments distincts du groupe $\drl_{U \supseteq Y^{(1)}} H^1(U_s, \Z / m\Z)$  avec $M > m^n*card(H^1(\CY_s, \mu_m))$. Pour chaque $\alpha_i$, il existe un ouvert $(U_i)_s$ de $\CY_s$ contenant $Y^{(1)}$ tel que $\alpha_i \in H^1((U_s)_i, \mu_m)$. Cependant,  comme $ card(H^1(\bigcap_{i=1}^M (U_s)_i, \mu_m)) \leq m^n*card(H^1(\CY_s, \mu_m)) $, il existe $\alpha_i, \alpha_j$ avec $i \neq j$ tels que leurs restrictions dans le groupe $ H^1(\bigcap_{i=1}^M (U_s)_i, \mu_m)$ sont égales par le principe des tiroirs, d'où une contradiction.\\

En vertu du diagramme commutatif suivant,
\[\begin{tikzcd}
  0 & H^1(l, \Z/m\Z) & H^1(Y, \Z/m\Z)  & \displaystyle{\drl_{U \subseteq Y^{(1)}}} H^1(U_s, \Z/m\Z) \\
  & \prod_{c \in Y^{(1)}} H^1(l(c), \Z /m\Z) &   \prod_{c \in Y^{(1)}} H^1(l(c), \Z /m\Z). \\
    \arrow[from=1-1, to=1-2]
    \arrow[from=1-2, to=1-3]
    \arrow[from=1-3, to=1-4]
    \arrow[from=1-2, to=2-2]
    \arrow[from=1-3, to=2-3]
    \arrow[from=2-2, to=2-3, "\simeq"]
\end{tikzcd}\]
on est donc ramené à prouver la finitude des groupes
$$
\CK^1(l, \Z / m\Z) \coloneqq Ker\bigg(H^1(l,\Z / m\Z) \ra \prod_{c \in Y^{(1)}} H^1(l(c), \Z /m\Z) \bigg)
$$
mais ces groupes sont égaux à $\Sha^1_n(l,\Z /m\Z)$ par la proposition \ref{prop1} et ces groupes sont finis par le lemme \ref{lemfinitude}, quitte à remplacer $l$ par une extension finie suffisamment grande de sorte que les composantes irréductibles de la fibre spéciale $Y_0$ soient géométriquement irréductibles afin que l'hypothèse \ref{hyp1}
soit trivialement satisfaite.\\

\section{Finitude du groupe $\Sha^1(K,G)$}\label{secfinitude}
Cette section a pour objet de faire le lien entre les différents résultats obtenus dans les parties précédentes. \\
Commençons par établir un petit lemme.
\begin{lemme}\label{lemcofinal}
Le morphisme canonique
$$
\drl_{V \supseteq X^{(1)}} H^1(V_l,G) \ra \drl_{U \supseteq Y^{(1)}} H^1(U,G)
$$
est un isomorphisme.
\end{lemme}

\proof Soit $U$ un ouvert de $\CY$ contenant $Y^{(1)}$ que l'on peut, quitte à le rétrécir, voir comme un ouvert de $Spec(R_l)$. Le morphisme $\varphi : \CY \ra \CX$ étant propre, l'ensemble $F \coloneqq \varphi(\CY \setminus U)$ est un fermé de $\CX$. Soit $V$ l'ouvert $\CX \setminus F$ de $\CX$, alors le schéma $V_l \coloneqq V \times_k l$ est un ouvert de $\CY$ inclus dans $U$. Il reste à vérifier que $V \supseteq X^{(1)}$. L'extension $l/k$ étant finie, le morphisme $\varphi$ est étale et fini : il préserve donc la codimension des points (cf \cite{Liu} proposition 3.23, chapitre 4). De fait, $F$ ne contient aucun point de $X^{(1)}$ par définition de $U$; cela implique que $V \supseteq X^{(1)}$ et à fortiori que $V_l \supseteq Y^{(1)}$. Cela nous permet de définir un morphisme naturel $\drl_{V \supseteq X^{(1)}} H^1(V_l,G) \ra \drl_{U \supseteq Y^{(1)}} H^1(U,G)$
et nous avons montré que l'ensemble des ouverts de la forme $V_l$ où $V$ est un ouvert de $\CX$ forme un système final et donc cofinal quand on considère la cohomologie, ce qui suffit.\\
\qed

Posons $\Gamma \coloneqq Gal(l/k)$, on dispose, en utilisant la suite spectrale d'Hochschild–Serre, pour tout ouvert $U \supseteq X^{(1)}$ de $\CX$, d'une suite exacte
$$
0 \ra H^1(\Gamma, G(U_{l})) \ra H^1(U, G) \ra H^1(U_{l}, G).
$$

Passer à la limite directe donne, en utilisant le lemme \ref{lemcofinal}, le diagramme commutatif suivant dont les suites verticales et horizontales sont exactes :

\[\begin{tikzcd}
    & 0 & 0 \\
	& \CK^1(X, G) & \CK^1(Y, G) \\
	{ \drl_{U \supseteq X^{(1)}} H^1(\Gamma, G(U_l))} & {H^1(X, G)} & {H^1(Y, G)} & {} \\
	 & {\prod_{c \in X^{(1)}} H^1(k(c), G)} & {\prod_{c \in Y^{(1)}} H^1(k(c), G).}  \\
	\arrow[from=2-2, to=2-3]
    \arrow[from=3-1, to=3-2]
    \arrow[from=3-2, to=3-3]
    \arrow[from=4-2, to=4-3]
	\arrow[from=2-2, to=3-2]
	\arrow[from=2-3, to=3-3]
    \arrow[from=1-2, to=2-2]
	\arrow[from=1-3, to=2-3]
    \arrow[from=3-2, to=4-2]
    \arrow[from=3-3, to=4-3]
\end{tikzcd}\]

Par le théorème d'Hilbert 90, le groupe $\CK^1(Y, G)$ est isomorphe à une produit fini de copies de $\Pic(Y)$ et de groupes de la forme $\CK^1(Y, \Z /m\Z)$. Ainsi, par les propositions \ref{propPicfini} et \ref{propdéployéfinitude}, le groupe $\ol{\CK^1(Y,G)}$ est de type fini et on entend montrer que le groupe $\ol{\CK^1(X,G)}$ l'est également.
 Une chasse au diagramme permet de montrer que le noyau du morphisme $\CK^1(X,G) \ra \CK^1(Y,G)$ est le groupe 
 $$
 D_{l} \coloneqq \text{Im}\bigg(\drl_{U \supseteq X^{(1)}} H^1(\Gamma, G(U_l)) \ra H^1(X, G) \bigg) \cap \CK^1(X, G).$$
 
Notons $D_{\CX}$ (respectivement $D_{\CY}$) le sous groupe divisible maximal de $\CK^1(X,G)$ (respectivement $\CK^1(Y,G)$).
Considérons  ensuite le diagramme commutatif suivant dont les suites verticales et horizontales sont exactes :

\[\begin{tikzcd}
   &  & 0 & 0 \\
  & &   D_{\CX} & D_{\CY} \\
	0 & D_{l} & \CK^1(X, G) & \CK^1(Y, G) \\
 & & {\ol{\CK^1(X, G)}} & {\ol{\CK^1(Y, G)}} \\
 &  & 0 & 0. \\
    \arrow[from=2-3, to=2-4, "f"]
    \arrow[from=3-1, to=3-2]
    \arrow[from=3-2, to=3-3]
    \arrow[from=3-3, to=3-4, "g"]
    \arrow[from=4-3, to=4-4, "h"]
	\arrow[from=2-3, to=3-3]
    \arrow[from=2-4, to=3-4]
    \arrow[from=3-1, to=3-2]
	\arrow[from=1-3, to=2-3]
    \arrow[from=1-4, to=2-4]
    \arrow[from=3-3, to=4-3]
    \arrow[from=3-4, to=4-4]
	\arrow[from=4-3, to=5-3]
    \arrow[from=4-4, to=5-4]
\end{tikzcd}\]

\begin{prop}\label{propDk1}
Supposons que l'hypothèse \ref{hyp1} soit satisfaite. Alors, le groupe $D_{l}$ est fini. 
\end{prop} 

\proof Tout d'abord, considérons le groupe $G(U_l)/ G(l)$ pour un ouvert $U \supseteq X^{(1)}$; il est isomorphe à un produit fini de $\G_m(U_l)/ \G_m(l)$ et de $(\Z/n\Z)(U_l)/ (\Z/n\Z)(l)$. \\
Par le lemme \ref{lemzariski}, les sections globales de $\CY$ et $U$ coïncident. Notons $\fm_l$ l'idéal maximal de $R_l$, alors
$$
\G_m(U_l)/ \G_m(l) \simeq R_l^*/l^* \simeq (l^* \times U^{(1)}) / l^* \simeq  U^{(1)}
$$où $ U^{(1)}$ désigne le groupe des unités principales, i.e.  $U^{(1)} = 1+\fm_l$. Ce dernier groupe est uniquement divisible par le lemme de Hensel, il en résulte que le groupe $H^1(\Gamma, \G_m(U_l)/ \G_m(l))$ est lui aussi uniquement divisible en plus d'être de torsion, il est donc trivial. Concernant les groupes de la forme $(\Z/n\Z)(U_l)/ (\Z/n\Z)(l)$, ils sont triviaux par connexité de $U_l$.\\

Ainsi, par la suite exacte
$$
0 \ra G(l) \ra G(U_l) \ra G(U_l)/ G(l) \ra 0
$$
 et au vu des considérations précédentes, on a la suite exacte cohomologique
$$
H^1(\Gamma, G(l)) \ra H^1(\Gamma,G(U_l)) \ra 0
$$
qui, après passage à la limite, devient 
$$
\drl_{U \supseteq X^{(1)}} H^1(\Gamma, G(l)) \ra \drl_{U \supseteq X^{(1)}} H^1(\Gamma,G(U_l)) \ra 0.
$$
On se ramène donc à établir la finitude du groupe  $\displaystyle \text{Im}\bigg(H^1(\Gamma, G(l)) \ra H^1(X, G)\bigg) \cap \CK^1(X, G) $.\\

Ce dernier groupe se plonge dans l'image dans $H^1(X, G)$ du groupe 
$$
\text{Ker}\bigg(H^1(k, G) \ra \prod_{c \in X^{(1)}} H^1(k(c), G)\bigg)
$$
qui par la proposition \ref{prop1} est égal au groupe
$$
\text{Ker}\bigg(H^1(k, G) \ra \prod_{c \in \CX^{(n)}} H^1(k(c), G)\bigg).
$$
\\
On en déduit la finitude du groupe $D_{l}$ dès lors que le groupe 
$\text{Ker}(H^1(k,G) \ra \prod_{i=1}^{n}H^1(k_i,G))$ est fini en vertu du lemme \ref{lemfinitude}. \qed\\

\begin{thme}\label{thme}
    Supposons que les hypothèses \ref{hyp1} et \ref{hyp2} sont satisfaites. Alors, le groupe $\Sha^1(K,G)$ est fini.
\end{thme}

\proof La proposition \ref{propDk1} donne la finitude du groupe $D_{l}$. D'une part, le groupe \text{Coker}(f) est de $m$-torsion finie pour tout entier $m$ : en effet, le groupe $D_{\CY}$ est de $m$-torsion finie d'après la proposition \ref{propPicfini}. Comme l'image d'un groupe divisible est également divisible, la suite exacte
$$
0 \ra Im(D_{\CX}) \ra D_{\CY} \ra \text{Coker}(f) \ra 0
$$
est scindée ce qui implique que le groupe $\text{Coker}(f)$ est isomorphe à un sous-groupe de $D_{\CY}$; sa $m$-torsion est donc bien finie.\\
D'autre part, on dispose par le lemme du serpent d'une suite exacte
$$
D_{l} \ra \text{Ker}(h) \ra \text{Coker}(f).
$$
qui démontre que la $m$-torsion du groupe $\text{Ker}(h)$ est finie pour tout entier $m$.\\

En outre, si l'on compose le morphisme de "restriction" (qui correspond à la limite des morphismes de restriction usuels) $h : \ol{\CK^1(X, G)} \ra \ol{\CK^1(Y, G)}$ par le morphisme de "corestriction" $\ol{\CK^1(Y, G)} \ra \ol{\CK^1(X, G)}$, on obtient la multiplication par $[l:k]$ sur $\ol{\CK^1(X, G)}$ (\cite{Milne},chap. \textit{V},§ 1, lemme 1.12). Il s'ensuit que le noyau de $h$ est d'exposant fini et est donc fini par ce qui précède; on en déduit au regard du fait que $\ol{\CK^1(Y, G)}$ est de type fini que le groupe $\ol{\CK^1(X,G)}$ l'est aussi.\\
Rappelons que $D_{\CX}$ désigne le sous-groupe divisible maximal de $\CK^1(X,G)$. Étant donné l'isomorphisme 
$$\CK^1(X,G) \simeq D_{\CX} \bigoplus \ol{\CK^1(X,G)}$$ 
et que l'image d'un groupe divisible dans un groupe d'exposant fini est nécessairement nulle, on peut alors conclure que l'image de $\CK^1(X,G)$ dans $H^1(K,G)$ est finie et donc que le groupe $\Sha^1(K,G)$ est fini par le lemme \ref{lemfinal}. \qed

\rmke L'hypothèse \ref{hyp2} ne saurait être nécessaire pour s'assurer de la véracité du théorème \ref{thme}. En effet, ce dernier est trivialement vrai lorsque $G = \gm$ et ce quelle que soit la géométrie de la fibre spéciale.\\

Concluons cette section avec un résultat qui nous permet de s'assurer que les hypothèses \ref{hyp1} et \ref{hyp2} sont raisonnables.
\begin{prop}
    Lorsque $R$ est de dimension $2$, les hypothèses \ref{hyp1} et \ref{hyp2} sont indépendantes de la désingularisation choisie.
\end{prop}

\proof Soient $\CX, \CX'$ deux désingularisations de $Spec(R)$. Par le théorème faible de factorisation pour les schémas quasi-excellents réguliers sur $\Q$ (\cite{AbTe19} 1.6), l'application birationnelle naturelle $\CX_1 \dra \CX_2$ se factorise en une composition d'éclatements et de contractions; il suffit donc d'établir que les hypothèses \ref{hyp1} et \ref{hyp2} sont préservées par un éclatement en un point fermé. Supposons donc que le morphisme $\CX' \ra \CX$ est un éclatement d'un point fermé $\beta \in \CX$.
\begin{enumerate}[label= \roman*)]

    \item Soient $\CX_k$ et $\CX'_k$ les fibres spéciales respectives de $\CX$ et $\CX'$ munies de leur structure réduite;  Notons $\CX_1,...,\CX_n$ (respectivement $\CX_1',...,\CX_n', \CX_{n+1}'$), les composantes irréductibles de $\CX_k$ (respectivement de $\CX_k')$. Puisque $\CX$ est régulier, on a que $\CX_i' = \CX_i$ pour tout $1 \leq i \leq n$ et $\CX_{n+1}' \simeq \P^1_{k(\beta)}$.
    Soit $P$ (respectivement $P'$) l'ensemble des points d'intersections des composantes irréductibles de $\CX_k$ (respectivement de $\CX'_k)$ et $i_p$ (respectivement $i_{p'})$ l'inclusion $p \in P\hookrightarrow \CX$ (respectivement l'inclusion $p' \in P'\hookrightarrow \CX')$. 

Les suites exactes 
$$
    0 \ra \CO_{\CX}^* \ra \prod_{i=1}^n \CO_{\CX_i}^* \ra \bigoplus_{p \in P} (i_p)_*\gm \ra 0
$$
et
$$
    0 \ra \CO_{\CX'}* \ra \prod_{i=1}^n \CO_{\CX'_i}^* \ra \bigoplus_{p' \in P'} (i_{p'})_*\gm \ra 0
$$

induisent les suites exactes suivantes 

\begin{equation}\label{eqsuite1}
\prod_{i=1}^n H^0(\CX_i, \gm) \ra \prod_{p \in P} H^0(k(p), \gm) \ra \Pic(\CX_k) \ra \prod_{i=1}^n \Pic(\CX_i)
\end{equation}

\begin{equation}\label{eqsuite2}
 \prod_{i=1}^{n+1} H^0(\CX'_i, \gm) \ra \prod_{p' \in P'} H^0(k(p), \gm) \ra \Pic(\CX'_k) \ra \prod_{i=1}^n \Pic(\CX'_i).
\end{equation}

On distingue deux situations différentes :

\begin{itemize}
    \item Supposons que $\beta$ est un point régulier de $\CX_k$; on peut sans perdre de généralité supposer que $\beta \in\CX_1$. Notons $p'_1$ le point d'intersection de $\CX_1$ et $\CX'_{n+1}$ dans $\CX'$ de sorte que $P' = P \cup \{p'_1\}$.\\
On définit un morphisme $$\prod_{i=1}^{n+1} H^0(\CX'_i, \gm) \ra \prod_{i=1}^n H^0(\CX_i, \gm)$$  par $(x_1,...,x_{n+1}) \mapsto (x_1,...,x_n)$ et un morphisme
$$\prod_{p' \in P'} H^0(k(p'), \gm) = \prod_{p' \in P' \setminus\{p'_1\}} H^0(k(p'), \gm) \times k(\beta)^* \ra \prod_{p \in P} H^0(k(p), \gm)$$
par $(\alpha_1,...,\alpha_l,\alpha_{p'_1}) \mapsto (\alpha_1,...,\alpha_l,\alpha_{p'_1})$.

Ces morphismes induisent un diagramme commutatif 
\[\begin{tikzcd}\label{diagrammeconoyau}
\prod_{i=1}^n H^0(\CX_i, \gm) & \prod_{p \in P} H^0(k(p), \gm)\\
 \prod_{i=1}^n H^0(\CX'_i, \gm) & \prod_{p' \in P'} H^0(k(p), \gm)
    \arrow[from=1-1, to=1-2]
    \arrow[from=2-2, to=1-2]
    \arrow[from=2-1, to=2-2]
    \arrow[from=2-1, to=1-1]
\end{tikzcd}\]
dont les conoyaux des morphismes horizontaux sont des isomorphismes. Il en résulte l'isomorphisme
\begin{equation}\label{equationPicard}
    Ker(\Pic(\CX_k) \ra \prod_{i=1}^n \Pic(\CX_i)) \simeq  Ker(\Pic(\CX'_k) \ra \prod_{i=1}^{n+1} \Pic(\CX_i')).
\end{equation}
 en vertu des suites exactes \ref{eqsuite1} et \ref{eqsuite2}.

\item Si $\beta$ est un point singulier de $\CX_k$, alors $\beta$ est un point d'intersection entre deux composantes irréductibles $\CX_1$ et $\CX_2$ de $\CX_k$. 
Notons $p'_i$ le point d'intersection entre $\CX_i$ et $\CX_{n+1}$ pour $i \in \{1,2\}$.
On définit un morphisme $$\prod_{i=1}^{n+1} H^0(\CX'_i, \gm) \ra \prod_{i=1}^n H^0(\CX_i, \gm)$$  par $(x_1,...,x_{n+1}) \mapsto (x_1,...,x_n)$ et un morphisme
$$\prod_{p' \in P'} H^0(k(p'), \gm) = \prod_{p' \in P' \setminus\{p'_1,p'_2\}} H^0(k(p'), \gm) \times k(\beta)^* \times k(\beta)^* \ra \prod_{p \in P \setminus \{\beta \}} H^0(k(p), \gm) \times k(\beta)^*$$
par $(\alpha_1,...,\alpha_l,\alpha_{p_1}, \alpha_{p_2}) \mapsto (\alpha_1,...,\alpha_l,\alpha_{p_1}\alpha_{p_2}^{-1})$.

Ces morphismes induisent eux aussi le diagramme commutatif \[\begin{tikzcd}\label{diagrammeconoyau}
\prod_{i=1}^n H^0(\CX_i, \gm) & \prod_{p \in P} H^0(k(p), \gm)\\
 \prod_{i=1}^n H^0(\CX'_i, \gm) & \prod_{p' \in P'} H^0(k(p), \gm)
    \arrow[from=1-1, to=1-2]
    \arrow[from=2-2, to=1-2]
    \arrow[from=2-1, to=2-2]
    \arrow[from=2-1, to=1-1]
\end{tikzcd}\]
dont les conoyaux des morphismes horizontaux sont des isomorphismes; il en résulte à nouveau l'isomorphisme \ref{equationPicard}.\\
Ainsi, la véracité de l'hypothèse \ref{hyp2} ne dépend pas du modèle régulier choisi.
\end{itemize}
    \item On conserve les notations du premier paragraphe. Notons $\CX_j$ l'une des composantes irréductibles à laquelle $\beta$ appartient. Alors, $k_j \coloneqq k(\CX_j) \cap k^s$ s'injecte dans $k(\beta)$ ce qui donne l'égalité $$Ker(H^1(k,G) \ra \prod_i^nH^1(k_i,G)) = Ker(H^1(k,G) \ra \prod_i^{n+1}H^1(k_i,G)).$$ 
L'hypothèse \ref{hyp1} ne dépend donc également pas du modèle régulier choisi.
\end{enumerate}
\qed

\section{Trivialité dans le cas régulier}
Ce résultat de finitude établi, on est amené à se poser la question suivante : les groupes de Tate-Shafarevich $\Sha^1(K,G)$ sont-ils nuls ? Le théorème suivant répond à cette question par l'affirmative lorsque l'anneau des entiers de $K$ est régulier.

\thme\label{thmezero} Supposons que l'anneau $R$ soit régulier. Alors $\Sha^1(K,G) = 0$.

\proof L'hypothèse de régularité sur $R$ implique que 
$$
\Sha_{n}^1(K,G) \coloneqq Ker \bigg( H^1(k, G) \ra \prod_{c \in \CX^{(n)}} H^1(k(c), G) \bigg) = 0
$$
car la fibre spéciale ne consiste que du point fermé de $Spec(R)$ dont le corps résiduel est $k$.\\
Par la proposition \ref{prop1}, le groupe 
$$
\Sha_{1}^1(G) \coloneqq Ker \bigg(H^1(k, G) \ra \prod_{c \in X^{(1)}} H^1(k(c), G) \bigg)
$$
est lui aussi trivial et, par la fin de la preuve de la proposition \ref{propDk1}, le groupe 
$$
D_l = \text{Im}(H^1(l, G(Y)) \ra H^1(X, G) ) \cap \CK^1(X, G)
$$
l'est également; il en résulte que le morphisme $g : \CK^1(X, G) \ra \CK^1(Y, G)$ est injectif puisque $D_{l} = Ker(g)$.\\

Le lemme suivant permet de terminer la preuve.

\begin{lemme}
Le groupe $\CK^1(Y, G)$ est nul.
\end{lemme}

 \proof Rappelons que le groupe $\CK^1(Y, G)$ est isomorphe à un produit fini de copies de $\Pic(Y)$ et de $\CK^1(Y, \Z / m\Z) \coloneqq Ker\bigg (H^1(Y, \Z/m\Z) \ra \prod_{c \in Y^{(1)}} H^1(l(c), \Z /m\Z) \bigg).$
\begin{itemize}
    \item Premièrement, à nouveau par régularité de $R_l$, le groupe $\Pic(U)$ est trivial pour tout ouvert $U$ de $Spec(R_l)$ d'où la trivialité de $\Pic(Y)$.
    \item Deuxièmement, soit $U$ un ouvert de $\CY$ contenant tous les points de codimension $1$ de $R_l$ et considérons le diagramme commutatif

\[\begin{tikzcd}
  0 & H^1(l, \Z/m\Z) & H^1(U, \Z/m\Z)  & H^1(U_s, \Z/m\Z) \\
  & \prod_{c \in Y^{(1)}} H^1(l(c), \Z /m\Z) &   \prod_{c \in Y^{(1)}} H^1(l(c), \Z /m\Z). \\
    \arrow[from=1-1, to=1-2]
    \arrow[from=1-2, to=1-3]
    \arrow[from=1-3, to=1-4]
    \arrow[from=1-2, to=2-2]
    \arrow[from=1-3, to=2-3]
    \arrow[from=2-2, to=2-3, "\simeq"]
\end{tikzcd}\]

La proposition \ref{prop1} montre que la première flèche verticale est injective. Concernant le groupe $H^1(U_s, \Z/m\Z) \simeq H^1(U_s, \mu_m)$, la suite de Kummer couplée à la nullité de $Pic(U_s)$ donne un isomorphisme $H^1(U_s, \mu_m) \simeq \CO_{U_s}(U_s)^*/m$. Quitte à rétrécir l'ensemble $U_s$, on peut le voir comme un ouvert de $R_s$ contenant tous les points de codimension $1$ de $R_s$. Par le lemme \ref{lemzariski}, les sections globales de $U_s$ coïncident avec celles de $R_s$ puisque $codim(R_s \setminus U_s) \geq 2$. En outre, l'anneau $R_s$ est hensélien ce qui implique que  $\CO_{U_s}(U_s)^*/m = 0$ et donc à fortiori que le groupe $H^1(U_s, \Z/m\Z)$ est nul. On conclut e passant à la limite directe que le groupe $\CK^1(\CY, \Z / m\Z)$ est nul.
\end{itemize} \qed

En conclusion, le groupe $\CK^1(Y, G)$ et à fortiori le groupe $\CK^1(X, G)$ sont triviaux. Il en découle que $\Sha^1(G)$, qui s'injecte dans l'image de $\CK^1(X, G)$ dans $H^1(K,G)$, est également nul. \qed \\

\section{Extensions à d'autres types de corps}
Au cours de la preuve du théorème \ref{thme}, le fait que le corps résiduel $k$ est un corps de nombre n'intervient qu'à deux reprises : dans la preuve du lemme \ref{lemfinitude} lorsque l'on utilise le lemme 2.2 de \cite{HarSza} qui affirme la finitude des groupes $\Sha^1_S(k,G)$ et lors de la preuve du lemme \ref{lemPicfinitude} où l'on fait appel au théorème de Mordell-Weil. Or, ce dernier est également valable pour des corps de type fini sur $\Q$ (voir par exemple le théorème 2.1 et l'exemple 2.4 de \cite{Conrad}); cela nous pousse à essayer d'étendre l'ensemble des corps étudiés à ceux dont le corps résiduel de l'anneau des entiers est de type fini sur $\Q$.

\subsection{Corps résiduel finiment engendré sur $\Q$}
Soit $k_0$ un corps de nombres. Considérons une $k_0$-algèbre $R$ locale, complète, normale et excellente de dimension $n$ dont le corps résiduel est un corps finiment engendré sur $\Q$ et soit $K$ son corps de fractions. Ainsi, $K$ est une extension finie de $k(Z)((x_1,...x_n))$ pour une $\Q$-variété $Z$ projective, lisse et géométriquement intègre. Concrètement, cela revient à regarder des corps similaires à ceux étudiés dans la première partie mais dont les corps résiduels de l'anneau des entiers sont de type fini sur $\Q$. Le seul résultat de la preuve du théorème \ref{thme} qui tombe en défaut est le lemme \ref{lemfinitude} qui permet de montrer que le groupe $\Sha^1_n(k,G)$ est fini. L'hypothèse suivante, qui est plus forte que l'hypothèse \ref{hyp1}, permet de contourner ce problème. 

\begin{hypoe}\label{hyp3}
On suppose que le modèle régulier $\CX \ra Spec(R)$ possède un $k(Z)$-point.
\end{hypoe}

En effet, si cette hypothèse est satisfaite, le groupe $\Sha^1_n(k(Z),G)$ est automatiquement trivial. Le reste de la preuve est identique en tout point et permet d'obtenir le théorème suivant :
 
\begin{thme}
    Supposons que les hypothèses \ref{hyp2} et \ref{hyp3} sont satisfaites. Alors, le groupe $\Sha^1(K,G)$ est fini.
\end{thme}

\subsection{Corps des fonctions sur des extensions finies de $\Q(\!(x)\!)$}
L'approche utilisée ici fonctionne également dans un autre contexte similaire : les corps des fonctions de variétés projectives lisses géométriquement intègres définies sur un corps $K = k_0(\!(x)\!)$ où $k_0$ désigne un corps de nombres. On note encore $R$ l'anneau $k_0[[x]]$ et $k = k_0$ le corps résiduel de $R$.\\

Soit $V$ une variété projective lisse géométriquement intègre de dimension $m$ définie sur $K$ et soit $L = K(V)$ son corps des fonctions. On considère un modèle propre $\CV$ sur $R$, i.e. tel que $\CV_{K} \simeq V$. De plus, on peut supposer que ce modèle est régulier par \cite{temkinDesingularizationQuasiexcellentSchemes2008} et que $\CV_k$ est un diviseur à croisement normaux strict de $\CV$. On travaillera à nouveau avec un groupe de type multiplicatif $G$ provenant d'un groupe de type multiplicatif $G_k$ défini sur $k$.

\rmke On peut se contenter de supposer que la $K$-variété $V$ est propre. En effet, si l'on considère un modèle $\CV_0$ de $V$ sur $Spec(R)$, le théorème de compactification de Nagata (\cite{Nagata}, §4, Theorem 2) garantit l'existence d'un schéma $\CV$ propre sur $Spec(R)$ tel que $\CV_0$ s'identifie à un ouvert dense de $\CV$; ainsi, la fibre générique de $\CV$ est bien isomorphe à $V$.\\

Le modèle régulier $\CV$ va jouer le rôle de la désingularisation $\CX \ra Spec(R)$ dans la première partie tandis que la fibre générique $V$ va jouer celui de l'ouvert régulier de $Spec(R)$. Étant donné la forte similarité avec le cadre de la première section, on se contente de donner les grandes lignes du raisonnement et d'expliciter les quelques points qui diffèrent.

\begin{lemme}
Pour tout point fermé $\alpha_{m+1}$ de $\CV$, il existe un modèle régulier $\CV'$ de $\CV$, un drapeau régulier $\beta = (\beta_{1},...,\beta_{m+1})$ dans $\CV'$ tel que $k(\alpha_{m+1}) \simeq k(\beta_{n+m})$ et une famille $\{S_2,...,S_{m+1}\}$ d'anneaux de valuation discrète henséliens tels que 
\begin{itemize}
    \item Pour tout $2\leq i < m+1$, on a  $k(\beta_i) \simeq Frac(S_{i+1})$.
    \item  Pour tout $2 < i \leq m+1$, on a $k(S_i) \simeq k(\beta_i)$ où $k(S_i)$ désigne le corps résiduel de $S_i$.
    \item $\beta_1 \in V^{(1)}$.
\end{itemize}
\end{lemme}

\proof La preuve est sensiblement identique à celle du lemme \ref{lemdimension}. Pour un point fermé $\beta_{m+1},$ on construit de la même façon un drapeau $\beta = \{\beta_{1},...,\beta_{m+1}\}$ dans un modèle régulier $\CV'$ tel que $\beta_m$ n'appartient pas à la fibre spéciale de $\CV'$; cela revient à dire que $\beta_m$ est envoyé sur le point générique de $Spec(R)$. Il en va alors nécessairement de même pour $\beta_1.$ \qed \\

Notons $W \coloneqq V^{(1)}$ l'ensemble des points de codimension $1$ de $\CV$ appartenant à la fibre générique $V$. On définit pour tout faisceau étale sur $\CV$
$$
H^i(W, \sF) \coloneqq \drl_{U \supseteq \CV^{(1)}} H^i(U,\sF)
$$
et
$$
\CK^1(W,G) \coloneqq Ker(H^1(W,G) \ra \prod_{c \in \CV^{(1)}} H^1(k(c), G)).
$$
On définit de manière analogue les groupes $H^i(W_l, \sF)$ et $\CK^1(W_l,G)$ pour une extension finie $l$ de $k$ suffisamment grande déployant $G$.

\begin{lemme}
    Le groupe $\Sha^1(L,G)$ s'injecte dans l'image de $\CK^1(W,G)$ dans $H^1(L,G)$.
\end{lemme}

\proof Idem que celle du lemme \ref{lemfinal}. \qed \\

On introduit à présent les analogue des hypothèses \ref{hyp1} \ref{hyp2} dans notre cadre. 

\begin{hypoe}\label{hyp5}
Soient $Z_1,...,Z_l$ les composantes irréductibles de la fibre spéciale $\CX_k$ et posons $k_i \coloneqq k(Z_i) \cap \overline{k}$. On suppose que le noyau du morphisme
$$
H^1(k, G) \ra \prod_{i=1}^l H^1(k_i, G)
$$
est fini.
\end{hypoe}

\begin{hypoe}\label{hyp4}
Soient $W_1,...,W_n$ les composantes irréductibles de la fibre spéciale $\CV_k$. On suppose que le noyau du morphisme $\Pic(\CV_k) \ra \prod_{i=1}^n\Pic(W_i)$ est fini.
\end{hypoe}

\rmke Les deux hypothèses \ref{hyp5} et
\ref{hyp4} sont vérifiées si la variété $V$ a bonne réduction en le point fermé de $Spec(R)$.\\

Comme $V$ est géométriquement intègre, on a que $H^0(V,\CO_V) = K$; le théorème de factorisation de Stein \ref{lemstein} montre alors que $H^0(\CV,\CO_{\CV}) \simeq R$. On peut alors suivre \textit{mutatis mutandis} le raisonnement employé dans la section \ref{secdéployé} pour montrer le lemme suivant.

\begin{lemme}
 Si l'hypothèse \ref{hyp5} est satisfaite, le groupe $\ol{\CK^1(W_l,G)}$ est de type fini et sa $m$-torsion est finie pour tout entier $m > 0$.
\end{lemme}

On peut enfin reprendre le raisonnement de la section \ref{secfinitude} pour obtenir l'analogue du théorème \ref{thme} :

\begin{thme}
Supposons que les hypothèses \ref{hyp5} et \ref{hyp4} sont satisfaites. Alors, le groupe $\Sha^1(L,G)$ est fini.
\end{thme}

\bibliographystyle{alpha}
\bibliography{Biblio}
\end{document}